\documentclass[12pt]{amsart}
\usepackage{amsmath}
\usepackage{amssymb}
\usepackage{amscd}
\usepackage[mathscr]{eucal}
\usepackage{latexsym} 
\usepackage[all]{xy}

\theoremstyle{definition}

\theoremstyle{remark}

\numberwithin{equation}{section}

\def\-rig{\text{\rm -rig}}
\def\-log{\text{\rm -log}}
\def\-dif{\text{\rm -dif}}

\def\lim{\text{\rm lim}}

\begin{document}
\begin{center}
\textbf{\large{Relations between Multi-Poly-Bernoulli numbers and Poly-Bernoulli numbers of negative index}}
\end{center}

\vspace{0.5cm}
\begin{center}
HIROYUKI $\hspace{3pt}$ KOMAKI
\end{center}

\vspace{0.5cm}

Abstract.
Poly-Bernoulli numbers $B_n^{(k)}\in\mathbb{Q}$\,($n \geq 0$,\,$k \in \mathbb{Z}$) are defined by Kaneko in 1997.
Multi-Poly-Bernoulli numbers\,$B_n^{(k_1,k_2,\ldots, k_r)}$, defined by using multiple polylogarithms, are
generations of Kaneko's Poly-Bernoulli numbers\,$B_n^{(k)}$. We researched relations between Multi-Poly-Bernoulli numbers and 
Poly-Bernoulli numbers of negative index in particular.
In section 2, we introduce a identity for Multi-Poly-Bernoulli numbers of negative index which was proved by Kamano.
In section 3, as main results, we introduce some relations between Multi-Poly-Bernoulli numbers and 
Poly-Bernoulli numbers of negative index in particular.

\section{Introduction}

For any integer $k$,\,Kaneko[1] introduced Poly-Bernoulli numbers of index $k$ by the following generating function:
\begin{center}
$\dfrac{Li_k(1-e^{-t})}{1-e^{-t}}$=$\displaystyle \sum_{n=0}^{\infty}{B_n^{(k)}}\dfrac{t^n}{n!}$,
\end{center}
where $Li_k(t)$ is the $k$-th polylogarithm defined by
\begin{center}
$Li_k(t)$=$\displaystyle\sum_{n=1}^{\infty}\dfrac{t^n}{n^k}$.
\end{center}
Since $Li_1(t)$=$-\mathrm{log}(1-t)$,\,the number $B_n^{(1)}$ is the ordinary $n$-th Bernoulli number $B_n$, which is defined by
\begin{center}
$\displaystyle\dfrac{te^t}{e^t-1}$=$\displaystyle\sum_{n=0}^{\infty}B_n\dfrac{t^n}{n!}$.
\end{center}
It is known that Poly-Bernoulli numbers of negative index are positive integers and we have a closed formula
\begin{center}
$B_n^{(-k)}$=$\displaystyle\sum_{j=0}^{min(n,k)}(j!)^2\genfrac{\{}{\}}{0pt}{}{n+1}{j+1}
\genfrac{\{}{\}}{0pt}{}{k+1}{j+1}$.
\end{center}
In particular,we have the following duality formula:
\vspace{5pt}
\begin{center}
$B_n^{(-k)}$=${B_k^{(-n)}}$ \hspace{0.2cm} ($k,\,n \geq 0$).
\end{center}
Moreover,these numbers have combinatorial applocations: see [2] and [3] for details.

As a generalization of Poly-Bernoulli numbers, Multi-Poly-Bernoulli numbers $B_n^{(k_1,k_2,\ldots, k_r)}$ are defined for integers $k_1, \ldots ,k_r$ by the generating function
\begin{center}
$\dfrac{Li_{k_1, \ldots ,k_r} (1-e^{-t})}{(1-e^{-t})^r}$=$\displaystyle \sum_{n=0}^{\infty}{B_n^{(k_1, \ldots ,k_r)}} \dfrac{t^n}{n!}$,
\end{center}
where $Li_{k_1,\ldots,k_r}(t)$ is a multiple polylogarithm defined by
\begin{center}
$Li_{k_1,\ldots,k_r}(t)$=$\displaystyle \sum_{0<m_1<\cdots<m_r}\dfrac{t^{m_r}}{m_1^{k_1}\cdots m_r^{k_r}}$.
\end{center}
When $r=1$, the number $B_n^{(k)}$ is Poly-Bernoulli numbers. When $r=1$ and $k_1=1$, the number $B_n^{(1)}$ is
the classical Bernoulli numbers. It is also known that we have the following duality formula [4] for Multi-Poly Bernoulli numbers:
\begin{center}
$B_n^{(0, \ldots, 0,-k)}$=${B_k^{(0, \ldots, 0,-n)}}$.
\end{center}

\section{Relations of Multi-Poly-Bernoulli numbers of negative index}

In this section, we introduce a identity for Multi-Poly-Bernoulli numbers of negative index which was proved by Kamano [4]. If  $(k_1,\ldots,k_r) \neq (0,\ldots,0)$, then Multi-Poly-Bernoulli numbers of negative index 
$B_n^{(-k_1,-k_2,\ldots, -k_r)}$ have the following expression.

\vspace{5pt}
\textbf{Theorem2.1.}
Let $r$ be a positive integer and let $k_1, \ldots, k_r$ be non-negative integers with $(k_1,\ldots,k_r) \neq (0,\ldots,0)$.
We put $k:=k_1+\cdots +k_r$. Then the following identity holds:
\begin{center}
$B_n^{(-k_1,\ldots,-k_r)}$=$\displaystyle \sum_{l=1}^{k} {\alpha}_l^{(k_1,\ldots,k_r)} (l+r)^n\cdots(A)$,
\end{center}
where ${\alpha}_l^{(k_1,\ldots,k_r)}$\,  $(1\le l \le k)$ are integers depending only on $k_1, \ldots, k_r$, and
they are inductively determined by the following recurrence relations:

(i)\,${\alpha}_l^{(k_1)}$=$(-1)^{l+k_1} l! \genfrac{\{}{\}}{0pt}{}{k_1}{l}$,

{\vspace{3pt}}

(ii)\,${\alpha}_l^{(k_1,\ldots,k_{r-1},0)}$=${\alpha}_l^{(k_1,\ldots,k_{r-1})}$,

{\vspace{3pt}}

(iii)\,${\alpha}_l^{(k_1,\ldots,k_{r-1},k_r+1)}$=$(l+r-1) {\alpha}_{l-1}^{(k_1,\ldots,k_r)}-l{\alpha}_l^{(k_1,\ldots,k_r)}$.

Here we set \begin{center}
${\alpha}_0^{(k_1,\ldots,k_r)}=\begin{cases}\
1 \hspace{5pt} \mathrm{if}\,(k_1,\ldots,k_r)=(0,\ldots,0),\\
{\hspace{3pt}}0 \hspace{5pt} \mathrm{otherwise},
\end{cases}$
\end{center} \hspace{5pt}  and
${\alpha}_l^{(k_1,\ldots,k_r)}$=$0$ {\hspace{3pt}} for $l>k$.

\hspace{-0.45cm} First we give recurrence relation [4] of Multi-Poly-Bernoulli numbers for the proof of the Theorem2.1.

\vspace{5pt}
\hspace{-0.45cm} \textbf{Lemma2.2.}
For non-negative integers $n,\,k_1,\ldots, k_r$, we have
\begin{center}
$B_n^{(-k_1,\ldots,-k_{r-1},-k_r-1)}$=$ \displaystyle \sum_{m=0}^{n} \genfrac{(}{)}{0pt}{}{n}{m}
B_{m+1}^{(-k_1,\ldots,-k_r)} +rB_n^{(-k_1,\ldots,-k_r)}-{B_{n+1}^{(-k_1,\ldots,-k_r)}}$.
\end{center}
Theorem2.1 is proved by induction on $r$. The following lemma [4] says that Theorem2.1 holds for $r=1$.

\vspace{5pt}
\hspace{-0.45cm} \textbf{Lemma2.3.}
For $n\ge{0}$ and $k\ge{1}$, we have
\begin{center}
$B_n^{(-k)}$=$\displaystyle \sum_{l=1}^{k} (-1)^{l+k} l! \genfrac{\{}{\}}{0pt}{}{k}{l}
{(l+1)^n}$.
\end{center}
We note here the Corollary2.4 has been proved by Hamahata and Masubuchi [5].

\vspace{5pt}
\hspace{-0.45cm}  \textbf{Corollary2.4.}
Let $r$ be a positive integer and let $n$ and $k$ be non-negative integers. Then the following identities hold:

(1) \hspace{3pt} ${B_n}^{\overbrace {\scriptstyle{(0, \ldots, 0 )}}^{r}}$=$r^n$,

\vspace{5pt}
(2) \hspace{3pt} (duality) \hspace{1pt}  $B_n^{(0, \ldots, 0,-k)}$=${B_k^{(0, \ldots, 0,-n)}}$,

\vspace{5pt}
(3) \hspace{3pt} $B_n^{(-k_1, \ldots, -k_{r-1},0)}$=$ \displaystyle \sum_{i=0}^{n}  {\binom{n}{i}} B_i^{(-k_1, \ldots, -k_{r-1})}$ $(r\ge{2})$,

\vspace{5pt}
(4) \hspace{3pt} $ \displaystyle \sum_{i=0}^{k} {\binom{k}{i}} B_n^{(-i,i-k)} p^i q^{k-i}$=
 $ \displaystyle \sum_{i=0}^{k}  \displaystyle \sum_{j=0}^{n}   {\binom{k}{i}}  {\binom{n}{j}} (p+q)^i q^{k-i} B_j^{(-i)} {B_{n-j}^{(i-k)}}$ where $p$

\hspace{0.45cm} and $q$ are any real numbers.

\vspace{5pt}
We use the following generating function of Multi-Poly-Bernoulli numbers of negative index for the proof of Corollary2.4.
This generating function is a natural generalization of the following function:
\begin{center}
$\displaystyle \sum_{n=0}^{\infty}  \displaystyle \sum_{k=0}^{\infty} B_n^{(-k)} \dfrac{x^n}{n!} \dfrac{y^k}{k!}=
\dfrac{1}{e^{-x}+e^{-y}-1}$.
\end{center}

\textbf{Theorem2.5.} The following identity holds:

\vspace{5pt}
\hspace{1cm} $\displaystyle \sum_{{k_1}=0}^{\infty} \cdots \displaystyle \sum_{{k_r}=0}^{\infty}  \displaystyle \sum_{n=0}^{\infty}
B_n^{(-k_1, \ldots, -k_r)} \dfrac{{x_1}^{k_1}}{{k_1}!} \cdots  \dfrac{{x_r}^{k_r}}{{k_r}!}  \dfrac{t^n}{n!}$

\vspace{5pt}
\hspace{0.5cm} =$\dfrac{1} {(e^{-x_1-x_2- \cdots -x_r}+e^{-t}-1)(e^{-x_2- \cdots -x_r}+e^{-t}-1) \cdots (e^{-x_r}+e^{-t}-1)}$.

\vspace{5pt}
We can express Multi-Poly-Bernoulli numbers and Poly-Bernoulli numbers of negative index in a sum of powers by using Theorem2.1.
We give examples [4], [5]  of Theorem2.1 for $1\le r \le 3$ and $1 \le k \le 3$ (Table 1).

{\vspace{5pt}}
\begin{center}
\begin{tabular}{c} \hline
$r=1$ \\ \hline 
$B_n^{(-1)}=2^n$ \\
$B_n^{(-2)}=-2^n+2\cdot 3^n$ \\
$B_n^{(-3)}=2^n-6 \cdot 3^n+6 \cdot 4^n$ \\ \hline
\end{tabular}
\end{center}
{\vspace{5pt}}
\begin{center}
\begin{tabular}{c} \hline
$r=2$ \\ \hline
$B_n^{(0,-1)}=2\cdot 3^n$\\
{\vspace{5pt}}
$B_n^{(-1,0)}=3^n$\\
$B_n^{(0,-2)}=-2\cdot 3^n+6\cdot 4^n$\\
{\vspace{5pt}}
$B_n^{(-2,0)}=-3^n+2\cdot 4^n$\\
$B_n^{(0,-3)}=2\cdot 3^n-18\cdot 4^n+24 \cdot 5^n$\\
$B_n^{(-1,-2)}=3^n-9\cdot 4^n+12\cdot 5^n$\\
$B_n^{(-2,-1)}=3^n-7\cdot 4^n+8\cdot 5^n$\\
$B_n^{(-3,0)}=3^n-6\cdot 4^n+6\cdot 5^n$\\ \hline
\end{tabular}
\end{center}
{\vspace{5pt}}
\begin{center}
\begin{tabular}{c} \hline
$r=3$ \\ \hline
$B_n^{(0,0,-1)}=3\cdot 4^n$\\
$B_n^{(0,-1,0)}=2\cdot 4^n$\\
{\vspace{5pt}}
$B_n^{(-1,0,0)}=4^n$\\
$B_n^{(0,0,-2)}=-3\cdot 4^n+12\cdot 5^n$\\
$B_n^{(0,-2,0)}=-2\cdot 4^n+6\cdot 5^n$\\
$B_n^{(-2,0,0)}=-4^n+2\cdot 5^n$\\
$B_n^{(0,-1,-1)}=-2\cdot 4^n+8\cdot 5^n$\\
$B_n^{(-1,0,-1)}=-4^n+4\cdot 5^n$\\
{\vspace{5pt}}
$B_n^{(-1,-1,0)}=-4^n+3\cdot 5^n$\\
$B_n^{(0,0,-3)}=3\cdot 4^n-36\cdot 5^n+60\cdot 6^n$\\
$B_n^{(0,-3,0)}=2\cdot 4^n-18\cdot 5^n+24\cdot 6^n$\\
$B_n^{(-3,0,0)}=4^n-6\cdot 5^n+6\cdot 6^n$\\
$B_n^{(0,-1,-2)}=2\cdot 4^n-24\cdot 5^n+40\cdot 6^n$\\
$B_n^{(0,-2,-1)}=2\cdot 4^n-20\cdot 5^n+30\cdot 6^n$\\
$B_n^{(-1,0,-2)}=4^n-12\cdot 5^n+20\cdot 6^n$\\
$B_n^{(-1,-2,0)}=4^n-9\cdot 5^n+12\cdot 6^n$\\
$B_n^{(-2,0,-1)}=4^n-8\cdot 5^n+10\cdot 6^n$\\
$B_n^{(-2,-1,0)}=4^n-7\cdot 5^n+8\cdot 6^n$\\
$B_n^{(-1,-1,-1)}=4^n-10\cdot 5^n+15\cdot 6^n$\\ \hline
\end{tabular}
\end{center}

We found regularities from Table 1 and got the following relations. The proof uses Theorem2.1.

\vspace{5pt}
\textbf{Theorem2.6.} We have the following relations, and i-th component is $-1$ and others are $0$ in (3).

(1)$B_n^{(-1,\overbrace {\scriptstyle{0,\ldots,0}}^{r-1})}$=$(r+1)^n$
=${B_n^{\overbrace {\scriptstyle{(0,\ldots,0)}}^{r+1}}}$.

(2)$B_n^{\overbrace {\scriptstyle{(0,\ldots,0}}^{r-1},-1)}$=${r(r+1)^n}$.

(3)$B_n^{\overbrace{\scriptstyle{(0,\ldots,0,-1,0,\ldots,0)}}^r}$=$i(r+1)^n \hspace{3pt} (1\le i \le r)$.

\vspace{5pt}
Proof. (1)We use Theorem2.1(A) and (i),(ii). The second equality obtains from Corollary1.4(1) .
\begin{align*}
B_n^{(-1,\overbrace {\scriptstyle{0,\ldots,0}}^{r-1})}
&=\displaystyle \sum_{l=1}^{1} {\alpha}_l^{(1,\overbrace {\scriptstyle{0,\ldots,0}}^{r-1})} (l+r)^n\\
&={\alpha}_1^{(1,\overbrace {\scriptstyle{0,\ldots,0}}^{r-1})}(1+r)^n\\
&={\alpha}_1^{(1)}(1+r)^n .
\end{align*}
Since ${\alpha}_1^{(1)}$=$(-1)^{1+1} 1! \genfrac{\{}{\}}{0pt}{}{1}{1}$=$1$, we obtain $B_n^{(-1,\overbrace {\scriptstyle{0,\ldots,0}}^{r-1})}$=$(r+1)^n$
=$B_n^{\overbrace {\scriptstyle{(0,\ldots,0)}}^{r+1}}$.

\vspace{5pt}
(2)We use Theorem2.1(A) and (iii).
\begin{align*}
B_n^{\overbrace {\scriptstyle{(0,\ldots,0}}^{r-1},-1)}
&=\displaystyle \sum_{l=1}^{1} {\alpha}_l^{\overbrace {\scriptstyle{(0,\ldots,0}}^{r-1},1)} (l+r)^n\\
&={\alpha}_1^{\overbrace {\scriptstyle{(0,\ldots,0}}^{r-1},1)} (1+r)^n.
\end{align*}
Here,
${\alpha}_1^{\overbrace {\scriptstyle{(0,\ldots,0}}^{r-1},1)}
=(1+r-1) {\alpha}_0^{\overbrace {\scriptstyle{(0,\ldots,0)}}^{r}} -1\cdot {\alpha}_1^{\overbrace {\scriptstyle{(0,\ldots,0)}}^{r}}$

$=r {\alpha}_0^{\overbrace {\scriptstyle{(0,\ldots,0)}}^{r}} - {\alpha}_1^{\overbrace {\scriptstyle{(0,\ldots,0)}}^{r}}=r$.

\hspace{8pt} Thus, we obtain $B_n^{\overbrace {\scriptstyle{(0,\ldots,0}}^{r-1},-1)}$=$r(r+1)^n$.

\vspace{5pt}
(3)We use Theorem2.1(A) and (ii).
\begin{align*}
B_n^{\overbrace{\scriptstyle{(0,\ldots,0,-1,0,\ldots,0)}}^r}
&=\displaystyle \sum_{l=1}^{1} {\alpha}_l^{\overbrace{\scriptstyle{(0,\ldots,0,-1,0,\ldots,0)}}^r}(l+r)^n\\
&={\alpha}_1^{\overbrace{\scriptstyle{(0,\ldots,0,-1,0,\ldots,0)}}^r} (r+1)^n\\
&={\alpha}_1^{\overbrace {\scriptstyle{(0,\ldots,0}}^{i-1},1)} (r+1)^n.
\end{align*}

\vspace{-10pt}
Here from (2), we have $B_n^{\overbrace {\scriptstyle{(0,\ldots,0}}^{i-1},-1)}$=$i(i+1)^n$. Moreover from Theorem2.1(A),
we have
\begin{align*}
\hspace{0.6cm} B_n^{\overbrace {\scriptstyle{(0,\ldots,0}}^{i-1},-1)}
&=\displaystyle \sum_{l=1}^{1} {\alpha}_l^{\overbrace {\scriptstyle{(0,\ldots,0}}^{i-1},1)} (l+i)^n\\
&={\alpha}_1^{\overbrace {\scriptstyle {(0,\ldots,0}}^{i-1},1)} (i+1)^n.
\end{align*}

\vspace{-10pt}
\hspace{0.45cm} Since ${\alpha}_1^{\overbrace {\scriptstyle{(0,\ldots,0}}^{i-1},1)}$=$i$, we obtain
$B_n^{\overbrace{\scriptstyle{(0,\ldots,0,-1,0,\ldots,0)}}^r}$=$i(r+1)^n $.

\hspace{0.45cm} We have $B_n^{\overbrace {\scriptstyle{(0,\ldots,0}}^{r-1},-1)}$=$r B_n^{(-1,\overbrace{\scriptstyle {0,\ldots,0}}^{r-1})}$
from (1) and (2).
\hspace{7cm} ${\Box}$

\section{Relations between Multi-Poly-Bernoulli numbers and Poly-Bernoulli numbers of negative index}

In this section, we introduce relations between Multi-Poly-Bernoulli numbers and Poly-Bernoulli numbers of negative index.
Poly-Bernoulli numbers of negative index can express by using the Stirling numbers of the second kind $\genfrac{\{}{\}}{0pt}{}{n}{m}$ and Multi-Poly-Bernoulli numbers. Here Stirling numbers of the second kind are the number of ways to divide a set of 
$n$ elements into $m$ nonempty sets.

\vspace{5pt}
\textbf{Theorem3.1.}
Poly-Bernoulli numbers $B_n^{(-k)}$ ($n\ge 0,\,k\ge 1$) can express as follows;

\vspace{-5pt}
(1)$B_n^{(-k)}$=$(-1)^{k-1} B_n^{(0,0)} + \displaystyle \sum_{r=1}^{k-1}
(-1)^{r+k+1} r! \genfrac{\{}{\}}{0pt}{}{k}{r+1} B_n^{\overbrace {\scriptstyle{(0,\ldots,0}}^{r},-1)}$${_{\textbf.}}$

(2)$B_n^{(-k)}$=$(-1)^{k-1} B_n^{(0,0)} + \displaystyle \sum_{r=1}^{k-1}
(-1)^{r+k+1} (r+1)! \genfrac{\{}{\}}{0pt}{}{k}{r+1} B_n^{(-1,\overbrace {\scriptstyle{0,\ldots,0}}^{r})}$${_{\textbf.}}$

(3)$B_n^{(-k)}$=$(-1)^{k-1} B_n^{(0,0)} + \displaystyle \sum_{r=i}^{k+i-2}
(-1)^{r-i+k} \dfrac{(r-i+2)!}{i} \genfrac{\{}{\}}{0pt}{}{k}{r-i+2} B_n^{\overbrace{\scriptstyle{(0,\ldots,0,-1,0,\ldots,0)}}^{r-i+2}}$${_{\textbf.}}$

\vspace{5pt}
\textbf{Example3.2.}
We give examples of Theorem3.1(1) and (2) for $1 \le k \le 4$.

(1)\,$B_n^{(-1)}$=$B_n^{(0,0)}$

{\hspace{13pt}} $B_n^{(-2)}$=$-B_n^{(0,0)}+B_n^{(0,-1)}$

{\hspace{13pt}} $B_n^{(-3)}$=$B_n^{(0,0)}-3B_n^{(0,-1)}+2B_n^{(0,0,-1)}$

{\hspace{13pt}} $B_n^{(-4)}$=$-B_n^{(0,0)}+7B_n^{(0,-1)}-12B_n^{(0,0,-1)}+6B_n^{(0,0,0,-1)}$\\

(2)\,$B_n^{(-1)}$=$B_n^{(0,0)}$

{\hspace{13pt}} $B_n^{(-2)}$=$-B_n^{(0,0)}+2B_n^{(-1,0)}$

{\hspace{13pt}} $B_n^{(-3)}$=$B_n^{(0,0)}-6B_n^{(-1,0)}+6B_n^{(-1,0,0)}$

{\hspace{13pt}} $B_n^{(-4)}$=$-B_n^{(0,0)}+14B_n^{(-1,0)}-36B_n^{(-1,0,0)}+24B_n^{(-1,0,0,0)}$

\vspace{5pt}
Proof of the Theorem3.1.

(1)From Lemma2.3, we have

\hspace{0.3cm} $B_n^{(-k)}$=$\displaystyle \sum_{l=1}^{k} (-1)^{l+k} l! \genfrac{\{}{\}}{0pt}{}{k}{l} {(l+1)^n}_{\textbf.}$

\hspace{0.3cm} Here by putting $l-1=r$, we obtain
\begin{align*}
\hspace{0.7cm} B_n^{(-k)}&=\displaystyle \sum_{r=0}^{k-1} (-1)^{r+k-1} (r+1)! \genfrac{\{}{\}}{0pt}{}{k}{r+1} (r+2)^n\\
&=(-1)^{k-1} 2^n +\displaystyle \sum_{r=1}^{k-1} (-1)^{r+k-1} (r+1)! \genfrac{\{}{\}}{0pt}{}{k}{r+1} (r+2)^n.\end{align*}
\hspace{0.5cm} Since $B_n^{(0,0)}$=$2^n$ from Corollary 2.4(1), we have
\begin{center}
$B_n^{(-k)}$=$(-1)^{k-1}B_n^{(0,0)}+\displaystyle \sum_{r=1}^{k-1} (-1)^{r+k-1} r! \genfrac{\{}{\}}{0pt}{}{k}{r+1} (r+1)(r+2)^n.$
\end{center}
\hspace{0.5cm} Here from Theorem2.6(2), we have $B_n^{\overbrace{\scriptstyle {(0,\cdots,0}}^{r-1},-1)}$=$r(r+1)^n$. Thus we have
\begin{center}
$B_n^{(-k)}$=$(-1)^{k-1}B_n^{(0,0)}+\displaystyle \sum_{r=1}^{k-1} (-1)^{r+k-1} r! \genfrac{\{}{\}}{0pt}{}{k}{r+1} B_n^{\overbrace {\scriptstyle{(0,\ldots,0}}^{r},-1)},$
\end{center}
\hspace{0.5cm} and we obtain the identity of (1).

(2) The proof of (2) uses the proof of (1). In the proof of (1), we have 
\begin{center}
$B_n^{(-k)}$=$(-1)^{k-1}B_n^{(0,0)}+\displaystyle \sum_{r=1}^{k-1} (-1)^{r+k-1} (r+1)! \genfrac{\{}{\}}{0pt}{}{k}{r+1}(r+2)^n.$
\end{center}
\hspace{0.5cm} Here from Theorem2.6(1), we have $B_n^{(-1,\overbrace {\scriptstyle{0,\ldots,0)}}^{r-1}}$=$(r+1)^n.$
Hence we have
\begin{center}
$B_n^{(-k)}$=$(-1)^{k-1}B_n^{(0,0)}+\displaystyle \sum_{r=1}^{k-1} (-1)^{r+k-1} (r+1)! \genfrac{\{}{\}}{0pt}{}{k}{r+1} B_n^{(-1,\overbrace {\scriptstyle{0,\ldots,0)}}^{r-1}},$
\end{center}
\hspace{0.45cm} and we obtain the identity of (2).

(3)From Lemma2.3, we have

\hspace{0.3cm} $B_n^{(-k)}$=$\displaystyle \sum_{l=1}^{k} (-1)^{l+k} l! \genfrac{\{}{\}}{0pt}{}{k}{l} (l+1)^n.$

\hspace{0.3cm} Here by putting $l-1=r-i+1$, we obtain
\begin{align*}
\hspace{0.7cm} B_n^{(-k)}&=\displaystyle \sum_{r=i-1}^{k+i-2} (-1)^{r-i+k} (r-i+2)! \genfrac{\{}{\}}{0pt}{}{k}{r-i+2} (r-i+3)^n\\
&=(-1)^{k-1} 2^n +\displaystyle \sum_{r=i}^{k+i-2} (-1)^{r-i+k} (r-i+2)! \genfrac{\{}{\}}{0pt}{}{k}{r-i+2} {\{(r-i+2)+1}\}^n\\
&=(-1)^{k-1} B_n^{(0,0)} +\displaystyle \sum_{r=i}^{k+i-2} (-1)^{r-i+k} \dfrac{(r-i+2)!}{i} \genfrac{\{}{\}}{0pt}{}{k}{r-i+2} i{\{(r-i+2)+1}\}^n.
\end{align*}
\hspace{0.6cm} Here from Theorem2.6(3), we have $B_n^{\overbrace{\scriptstyle{(0,\ldots,0,-1,0,\ldots,0)}}^r}$=$i(r+1)^n.$
Thus we have 
\begin{center}
$B_n^{(-k)}=(-1)^{k-1} B_n^{(0,0)} +\displaystyle \sum_{r=i}^{k+i-2} (-1)^{r-i+k} \dfrac{(r-i+2)!}{i} \genfrac{\{}{\}}{0pt}{}{k}{r-i+2} B_n^{\overbrace{\scriptstyle{(0,\ldots,0,-1,0,\ldots,0)}}^{r-i+2}},$
\end{center}
\hspace{0.7cm} and we obtain the identity of (3). \hspace{7cm} ${\Box}$

\vspace{5pt}
In the Thorem3.1, when we replace $r \to 2r$ and $i \to r+1$, we obtain the identity of (1). When we replace $i \to 1$, we obtain
the identity of (2). Hence the identity of (3) is generalization of (1) and (2).

Futhermore we can also express Poly-Bernoulli numbers of negative index using the Stirling numbers of the first kind
$\begin{bmatrix}
n\\m
\end{bmatrix}$
and Multi-Poly-Bernoulli numbers.
Here Stirling numbers of the first kind are the number of permutations of $n$ letters (elements of the symmetric group of degree $n$) that consist of $m$ disjoint cycles.

\vspace{6pt}
\textbf{Corollary3.3.} We have the following relations

(1)\,$B_n^{(-k)}$=$(-1)^{k-1} B_n^{(0,0)} + \displaystyle \sum_{r=1}^{k-1}
(-1)^{r+k+1} r! 
\begin{bmatrix}
-r-1\\-k
\end{bmatrix}
B_n^{\overbrace {\scriptstyle{(0,\ldots,0}}^{r},-1)}$.
 
(2)\,$B_n^{(-k)}$=$(-1)^{k-1} B_n^{(0,0)} + \displaystyle \sum_{r=1}^{k-1}
(-1)^{r+k+1} (r+1)! 
\begin{bmatrix}
-r-1\\-k
\end{bmatrix}
B_n^{(-1,\overbrace {\scriptstyle{0,\ldots,0}}^{r})}$.

(3)\,$B_n^{(-k)}$=$(-1)^{k-1} B_n^{(0,0)} + \displaystyle \sum_{r=i}^{k+i-2}
(-1)^{r-i+k} \dfrac{(r-i+2)!}{i} 
\begin{bmatrix}
-r+i-2\\-k
\end{bmatrix}
B_n^{\overbrace{\scriptstyle{(0,\ldots,0,-1,0,\ldots,0)}}^{r-i+2}}$.

\vspace{5pt} 
The proof of Corollary3.3 can be obtained from the following Lemma3.4 [1].

\vspace{5pt}
\textbf{Lemma3.4.} For any integers $n$ and $m$, we have
\begin{center}
$\begin{bmatrix}
n\\m
\end{bmatrix}
$=$\genfrac{\{}{\}}{0pt}{}{-m}{-n}.$
\end{center}

Next we see the sum of coefficients on Multi-Poly-Bernoulli numbers of the identity which hold on Theorem3.1.
Therefore we revisit Example3.2.

\vspace{5pt}
\textbf{Example3.2} (Example3.2 revisited). 

We give examples of Theorem3.1(1) and (2) for $1 \le k \le 4$.

(1)\,$B_n^{(-2)}$=$-B_n^{(0,0)}+B_n^{(0,-1)}$

{\hspace{10pt}} $B_n^{(-3)}$=$B_n^{(0,0)}-3B_n^{(0,-1)}+2B_n^{(0,0,-1)}$

{\hspace{10pt}} $B_n^{(-4)}$=$-B_n^{(0,0)}+7B_n^{(0,-1)}-12B_n^{(0,0,-1)}+6B_n^{(0,0,0,-1)}$\\

(2)\,$B_n^{(-2)}$=$-B_n^{(0,0)}+2B_n^{(-1,0)}$

{\hspace{10pt}} $B_n^{(-3)}$=$B_n^{(0,0)}-6B_n^{(-1,0)}+6B_n^{(-1,0,0)}$

{\hspace{10pt}} $B_n^{(-4)}$=$-B_n^{(0,0)}+14B_n^{(-1,0)}-36B_n^{(-1,0,0)}+24B_n^{(-1,0,0,0)}$

\vspace{5pt}
In the case of (1), the sum of coefficients on Multi-Poly-Bernoulli numbers are $0$ ($k \geq 2$).
In the case of (2), the sum of coefficients on Multi-Poly-Bernoulli numbers are $1$ ($k \geq 2$).
From here we can be considered the following relations.

\vspace{5pt}
\textbf{Theorem3.5.} We have the following relations for $k \geq 2$

(1)\,$(-1)^{k-1}+ \displaystyle \sum_{r=1}^{k-1} (-1)^{r+k+1} r! \genfrac{\{}{\}}{0pt}{}{k}{r+1}$=$0$.

(2)\,$(-1)^{k-1}+ \displaystyle \sum_{r=1}^{k-1} (-1)^{r+k+1} (r+1)! \genfrac{\{}{\}}{0pt}{}{k}{r+1}$=$1$.

(3)\,$(-1)^{k-1}+ \displaystyle \sum_{r=i}^{k+i-2} (-1)^{r-i+k} \dfrac{(r-i+2)!}{i}
\genfrac{\{}{\}}{0pt}{}{k}{r-i+2}$=$
\begin{cases}\
1 & \text{($k$: odd)}\\
\dfrac{2}{i} -1 & \text{($k$: even)}.
\end{cases}$

\vspace{5pt}
We regard the sums of coefficients as $1$ for $k=1$.

\vspace{5pt}
Proof. (1) We have
\begin{center}
$(-1)^{k-1}+\displaystyle \sum_{r=1}^{k-1} (-1)^{r+k+1} r! \genfrac{\{}{\}}{0pt}{}{k}{r+1}$=
$(-1)^{k-1}\biggl(1+\displaystyle \sum_{r=1}^{k-1} (-1)^r r! \genfrac{\{}{\}}{0pt}{}{k}{r+1}\biggr).$
\end{center}

\vspace{5pt}
Here, $(-1)^r r! \genfrac{\{}{\}}{0pt}{}{k}{r+1}$=$1$ for $r=0$ and $(-1)^r r! \genfrac{\{}{\}}{0pt}{}{k}{r+1}$=$0$ for $r=k$. 
Thus we have
\begin{align*}
\hspace{0.9cm} (-1)^{k-1}\biggl(1+\displaystyle \sum_{r=1}^{k-1} (-1)^r r! \genfrac{\{}{\}}{0pt}{}{k}{r+1}\biggr)
&=(-1)^{k-1}\biggl(1+\displaystyle \sum_{r=0}^{k} (-1)^r r! \genfrac{\{}{\}}{0pt}{}{k}{r+1}-1\biggr)\\
&=(-1)^{k-1}\displaystyle \sum_{r=0}^{k} (-1)^r r! \genfrac{\{}{\}}{0pt}{}{k}{r+1}\\
&=(-1)^{k-1}\displaystyle \sum_{r=0}^{k} (-1)^r
\begin{bmatrix}
r+1\\1
\end{bmatrix}
\genfrac{\{}{\}}{0pt}{}{k}{r+1}.
\end{align*}
Futhermore, since $\displaystyle \sum_{l=0}^{n} (-1)^l \genfrac{\{}{\}}{0pt}{}{n}{l}
\begin{bmatrix}
l\\m
\end{bmatrix}$=
$(-1)^m \delta_{m,n}$ ([1]), we have
\begin{align*}
(-1)^{k}\displaystyle \sum_{r=0}^{k} (-1)^{r+1}
\begin{bmatrix}
r+1\\1
\end{bmatrix}
\genfrac{\{}{\}}{0pt}{}{k}{r+1}
&=(-1)^k \delta_{1,k}\\
&=0,
\end{align*}
and we obtain the results.

(2) Considering in the same way with (1), we have
\begin{align*} (-1)^{k-1}+ \displaystyle \sum_{r=1}^{k-1} (-1)^{r+k+1} (r+1)! \genfrac{\{}{\}}{0pt}{}{k}{r+1}
&=(-1)^{k-1}\biggl(1+\displaystyle \sum_{r=1}^{k-1} (-1)^r (r+1)! \genfrac{\{}{\}}{0pt}{}{k}{r+1}\biggr)
\end{align*}
\begin{align*}
\hspace{7cm} &=(-1)^{k-1}\displaystyle \sum_{r=0}^{k-1} (-1)^r (r+1)! \genfrac{\{}{\}}{0pt}{}{k}{r+1}\\
&=(-1)^{k-1}\displaystyle \sum_{r=0}^{k-1} (-1)^r \displaystyle \sum_{l=0}^{r+1}
\begin{bmatrix}
r+1\\l
\end{bmatrix}
\genfrac{\{}{\}}{0pt}{}{k}{r+1}\\
&=(-1)^{k-1}\displaystyle \sum_{r=0}^{k}
(-1)^r\displaystyle \sum_{l=0}^{k}
\begin{bmatrix}
r+1\\l
\end{bmatrix} \genfrac{\{}{\}}{0pt}{}{k}{r+1} \\
&=(-1)^k \displaystyle \sum_{l=0}^{k} 
\displaystyle \sum_{r=0}^{k} (-1)^{r+1} \genfrac{\{}{\}}{0pt}{}{k}{r+1}
\begin{bmatrix}
r+1\\l 
\end{bmatrix}_{{\textbf.}}
\end{align*}
Here we use the aforesaid formula again; $\displaystyle \sum_{l=0}^{n} (-1)^l \genfrac{\{}{\}}{0pt}{}{n}{l}
\begin{bmatrix}
l\\m
\end{bmatrix}$=
$(-1)^m \delta_{m,n}$ ([1]). Then we have
\begin{align*}
(-1)^k \displaystyle \sum_{l=0}^{k} 
\displaystyle \sum_{r=0}^{k}(-1)^{r+1} \genfrac{\{}{\}}{0pt}{}{k}{r+1}
\begin{bmatrix}
r+1\\l
\end{bmatrix}
&=(-1)^k\displaystyle \sum_{l=0}^{k} (-1)^l \delta_{l,k}\\
&=(-1)^k \cdot (-1)^k\\
&=1.
\end{align*}

(3)\,(i) If $k$ is odd, we put $k=2m+1$. Then we have 
\begin{center}
$(-1)^{2m} +\displaystyle \sum_{r=i}^{2m+i-1} (-1)^{r-i+2m+1} \dfrac{(r-i+2)!}{i} 
\genfrac{\{}{\}}{0pt}{}{2m+1}{r-i+2}$
=$1+\displaystyle \sum_{r=i}^{2m+i-1} (-1)^{r-i+2m+1} \dfrac{(r-i+2)!}{i} 
\genfrac{\{}{\}}{0pt}{}{2m+1}{r-i+2}$.
\end{center}

Hence it suffices to show the following identity;
\begin{center}
$\displaystyle \sum_{r=i}^{2m+i-1} (-1)^{2m+r-i+1} \dfrac{(r-i+2)!}{i} 
\genfrac{\{}{\}}{0pt}{}{2m+1}{r-i+2}$=$0$.
\end{center}

\hspace{0.3cm} $\displaystyle \sum_{r=i}^{2m+i-1} (-1)^{2m+r-i+1} \dfrac{(r-i+2)!}{i} 
\genfrac{\{}{\}}{0pt}{}{2m+1}{r-i+2}$

$=\displaystyle \sum_{r=i-2}^{2m+i-1} (-1)^{2m+r-i+1} \dfrac{(r-i+2)!}{i} 
\genfrac{\{}{\}}{0pt}{}{2m+1}{r-i+2} -\dfrac{1}{i} \genfrac{\{}{\}}{0pt}{}{2m+1}{1}$

$=\dfrac{1}{i} \displaystyle \sum_{r=i-2}^{2m+i-1} (-1)^{2m+r-i+1} (r-i+2)! 
\genfrac{\{}{\}}{0pt}{}{2m+1}{r-i+2} -\dfrac{1}{i}$.

Here from Theorem3.5(2), since $(-1)^{k-1}+ \displaystyle \sum_{r=1}^{k-1} (-1)^{r+k+1} (r+1)! \genfrac{\{}{\}}{0pt}{}{k}{r+1}$=$1$, we put $k=2m+1$. Then we obtain
\begin{align*}
(-1)^{2m}+ \displaystyle \sum_{r=1}^{2m} (-1)^{r+2m} (r+1)! \genfrac{\{}{\}}{0pt}{}{2m+1}{r+1}&=1\\
\displaystyle \sum_{r=1}^{2m} (-1)^{r+2m} (r+1)! \genfrac{\{}{\}}{0pt}{}{2m+1}{r+1}&=0\\
\displaystyle \sum_{r=i}^{2m+i-1} (-1)^{2m+r-i+1} (r-i+2)! \genfrac{\{}{\}}{0pt}{}{2m+1}{r-i+2}&=0\\
\displaystyle \sum_{r=i-2}^{2m+i-1} (-1)^{2m+r-i+1} (r-i+2)! \genfrac{\{}{\}}{0pt}{}{2m+1}{r-i+2}-1&=0\\
\displaystyle \sum_{r=i-2}^{2m+i-1} (-1)^{2m+r-i+1} (r-i+2)! \genfrac{\{}{\}}{0pt}{}{2m+1}{r-i+2}&=1.
\end{align*}

Hence, we have $\displaystyle \sum_{r=i}^{2m+i-1} (-1)^{2m+r-i+1} \dfrac{(r-i+2)!}{i} 
\genfrac{\{}{\}}{0pt}{}{2m+1}{r-i+2}$=$0$ and, the sum of the coefficients are $1$.

(ii)\,If $k$ is even, we put $k=2m$. Then we have 

\begin{center}
$(-1)^{2m-1} +\displaystyle \sum_{r=i}^{2m+i-2} (-1)^{r-i+2m} \dfrac{(r-i+2)!}{i} 
\genfrac{\{}{\}}{0pt}{}{2m}{r-i+2}$
=$-1+\displaystyle \sum_{r=i}^{2m+i-2} (-1)^{r-i+2m} \dfrac{(r-i+2)!}{i} 
\genfrac{\{}{\}}{0pt}{}{2m}{r-i+2}$.
\end{center}

Hence, it suffices to show the following identity;
\begin{center}
$\displaystyle \sum_{r=i}^{2m+i-2} (-1)^{2m+r-i} \dfrac{(r-i+2)!}{i} 
\genfrac{\{}{\}}{0pt}{}{2m}{r-i+2}$=$\dfrac{2}{i}$.
\end{center}

\hspace{0.25cm} $\displaystyle \sum_{r=i}^{2m+i-2} (-1)^{2m+r-i} \dfrac{(r-i+2)!}{i} 
\genfrac{\{}{\}}{0pt}{}{2m}{r-i+2}$

$=\displaystyle \sum_{r=i-2}^{2m+i-2} (-1)^{2m+r-i} \dfrac{(r-i+2)!}{i} 
\genfrac{\{}{\}}{0pt}{}{2m}{r-i+2} +\dfrac{1}{i} \genfrac{\{}{\}}{0pt}{}{2m}{1}$

$=\dfrac{1}{i} \displaystyle \sum_{r=i-2}^{2m+i-2} (-1)^{2m+r-i} (r-i+2)! 
\genfrac{\{}{\}}{0pt}{}{2m}{r-i+2} +\dfrac{1}{i}$.

Here from Theorem3.5(2), we put $k=2m$. Then we have 
\begin{align*}
(-1)^{2m-1}+ \displaystyle \sum_{r=1}^{2m-1} (-1)^{r+2m+1} (r+1)! \genfrac{\{}{\}}{0pt}{}{2m}{r+1}&=1\\
\displaystyle \sum_{r=1}^{2m-1} (-1)^{r+2m+1} (r+1)! \genfrac{\{}{\}}{0pt}{}{2m}{r+1}&=2\\
\displaystyle \sum_{r=i}^{2m+i-2} (-1)^{2m+r-i} (r-i+2)! \genfrac{\{}{\}}{0pt}{}{2m}{r-i+2}&=2\\
\displaystyle \sum_{r=i-2}^{2m+i-2} (-1)^{2m+r-i} (r-i+2)! \genfrac{\{}{\}}{0pt}{}{2m}{r-i+2}+1&=2\\
\displaystyle \sum_{r=i-2}^{2m+i-2} (-1)^{2m+r-i} (r-i+2)! \genfrac{\{}{\}}{0pt}{}{2m}{r-i+2}&=1.
\end{align*}
Hence, since $\displaystyle \sum_{r=i}^{2m+i-2} (-1)^{2m+r-i} \dfrac{(r-i+2)!}{i} 
\genfrac{\{}{\}}{0pt}{}{2m}{r-i+2}$=$\dfrac{2}{i}$, the sum of coefficients are $\dfrac{2}{i}-1$. \hspace{12cm} ${\Box}$

\vspace{5pt}
We found that Poly-Bernoulli numbers of negative index can express using the Stirling numbers of the second kind $\genfrac{\{}{\}}{0pt}{}{n}{m}$ and the sum of Multi-Poly-Bernoulli numbers. 
This time, we introduce that special values of Multi-Poly-Bernoulli numbers which hold on Theorem2.6 can express
by using the sum of Poly-Bernoulli numbers.

\textbf{Theorem3.6}($r \geq 1$). We have the following relations

(1)\,$B_n^{(-1,\overbrace {\scriptstyle{0,\ldots,0}}^{r-1})}$=$ \dfrac{1}{r!} \displaystyle \sum_{k=1}^{r}
\begin{bmatrix}
r\\k
\end{bmatrix}
B_n^{(-k)}$.

(2)\,$B_n^{\overbrace {\scriptstyle{(0,\ldots,0}}^{r-1},-1)}$=$\dfrac{1}{(r-1)!} \displaystyle \sum_{k=1}^{r}
\begin{bmatrix}
r\\k
\end{bmatrix}
B_n^{(-k)}$.

(3)$B_n^{\overbrace{\scriptstyle{(0,\ldots,0,-1,0,\ldots,0)}}^r}$=$\dfrac{i}{r!}  \displaystyle \sum_{k=1}^{r}
\begin{bmatrix}
r\\k
\end{bmatrix}
B_n^{(-k)}$.

\vspace{5pt}
Proof of Theorem3.6.
\begin{align*}
(1) \dfrac{1}{r!} \displaystyle \sum_{k=1}^{r}
\begin{bmatrix}
r\\k
\end{bmatrix}
B_n^{(-k)}
&=\dfrac{1}{r!} \displaystyle \sum_{k=1}^{r}
\begin{bmatrix}
r\\k
\end{bmatrix}
\displaystyle \sum_{l=1}^{k} (-1)^{l+k} l! \genfrac{\{}{\}}{0pt}{}{k}{l} (l+1)^n\\
&=\dfrac{1}{r!} \displaystyle \sum_{k=0}^{r}
\begin{bmatrix}
r\\k
\end{bmatrix}
\displaystyle \sum_{l=0}^{k} (-1)^{l+k} l! \genfrac{\{}{\}}{0pt}{}{k}{l} (l+1)^n\\
&=\dfrac{1}{r!} \displaystyle \sum_{k=0}^{r}
\begin{bmatrix}
r\\k
\end{bmatrix}
\displaystyle \sum_{l=0}^{r} (-1)^{l+k} l! \genfrac{\{}{\}}{0pt}{}{k}{l} (l+1)^n\\
&=\dfrac{1}{r!} \displaystyle \sum_{l=0}^{r}(-1)^l l! (l+1)^n 
\displaystyle \sum_{k=l}^{r} (-1)^{k} 
\begin{bmatrix}
r\\k
\end{bmatrix}
\genfrac{\{}{\}}{0pt}{}{k}{l}\\
&=\dfrac{1}{r!} \displaystyle \sum_{l=0}^{r}(-1)^l l! (l+1)^n (-1)^l \delta_{l,r}\\
&=\dfrac{1}{r!} \displaystyle \sum_{l=0}^{r} l! (l+1)^n  \delta_{l,r}\\
&=\dfrac{1}{r!} \cdot r! (r+1)^n\\
&=(r+1)^n=B_n^{(-1,\overbrace {\scriptstyle{0,\ldots,0}}^{r-1})}.
\end{align*}
(2) We consider in the same way as (1), and we obtain
\begin{align*}
\dfrac{1}{(r-1)!} \displaystyle \sum_{k=1}^{r}
\begin{bmatrix}
r\\k
\end{bmatrix}B_n^{(-k)}&=r(r+1)^n=B_n^{\overbrace {\scriptstyle{(0,\ldots,0}}^{r-1},-1)}.
\end{align*}
(3) We consider in the same way as (1), and we obtain 
\begin{center}
$\dfrac{i}{r!} \displaystyle \sum_{k=1}^{r}
\begin{bmatrix}
r\\k
\end{bmatrix}B_n^{(-k)}=i(r+1)^n=B_n^{\overbrace{\scriptstyle{(0,\ldots,0,-1,0,\ldots,0)}}^r}$.
\end{center}
This completes the proof. \hspace{9cm} ${\Box}$

The following Collorary3.7 can be obtained by using Lemma3.4 in the identity of Theorem3.6.
Therefore we omit the proof.

\textbf{Collorary3.7}($ r\geq 1$). We have the following relations

(1)\,$B_n^{(-1,\overbrace{\scriptstyle {0,\ldots,0}}^{r-1})}$=$ \dfrac{1}{r!} \displaystyle \sum_{k=1}^{r}
\genfrac{\{}{\}}{0pt}{}{-k}{-r} B_n^{(-k)}$.

(2)\,$B_n^{\overbrace {\scriptstyle{(0,\ldots,0}}^{r-1},-1)}$=$\dfrac{1}{(r-1)!} \displaystyle \sum_{k=1}^{r}
\genfrac{\{}{\}}{0pt}{}{-k}{-r} B_n^{(-k)}$.

(3)$B_n^{\overbrace{\scriptstyle{(0,\ldots,0,-1,0,\ldots,0)}}^r}$=$\dfrac{i}{r!}  \displaystyle \sum_{k=1}^{r}
\genfrac{\{}{\}}{0pt}{}{-k}{-r} B_n^{(-k)}$.

\textbf{Example3.8.}

We give examples of Theorem3.6 for $2 \le r \le 4$.

\hspace{12pt}$B_n^{(-1,0)}$=$\dfrac{1}{2} B_n^{(-1)}+\dfrac{1}{2} B_n^{(-2)}$

\hspace{6pt}$B_n^{(-1,0,0)}$=$\dfrac{1}{3} B_n^{(-1)} +\dfrac{1}{2} B_n^{(-2)} +\dfrac{1}{6} B_n^{(-3)}$

$B_n^{(-1,0,0,0)}$=$\dfrac{1}{4} B_n^{(-1)} +\dfrac{11}{24} B_n^{(-2)} +\dfrac{1}{4} B_n^{(-3)}
+\dfrac{1}{24} B_n^{(-4)}$

\vspace{0.5cm}
\hspace{12pt}$B_n^{(0,-1)}$=$B_n^{(-1)}+B_n^{(-2)}$

\hspace{6pt}$B_n^{(0,0,-1)}$=$B_n^{(-1)} +\dfrac{3}{2}B_n^{(-2)} +\dfrac{1}{2} B_n^{(-3)}$

$B_n^{(0,0,0,-1)}$=$B_n^{(-1)}+\dfrac{11}{6} B_n^{(-2)} +B_n^{(-3)}+\dfrac{1}{6} B_n^{(-4)}$

Moreover, we consider $B_n^{(-m,\overbrace {\scriptstyle{0,\cdots,0}}^{r-1})}$ which are the generalizations of Theorem3.6(1).
First, we consider in the case of $m=2$, that is, $B_n^{(-2,\overbrace {\scriptstyle{0,\cdots,0}}^{r-1})}$.
We give examples for $1 \leq r \leq 4$ and $m=2$;
\begin{center}
\hspace{14pt} $B_n^{(-2)}$=$-2^n+2\cdot 3^n$

\hspace{7pt} $B_n^{(-2,0)}$=$-3^n+2\cdot4^n$

\hspace{3pt} $B_n^{(-2,0,0)}$=$-4^n+2\cdot5^n$

$B_n^{(-2,0,0,0)}$=$-5^n+2\cdot6^n$

$\cdot\cdot\cdot$.
\end{center}
We use Theorem2.6(1)\,($B_n^{(-1,\overbrace {\scriptstyle{0,\ldots,0}}^{r-1})}$=$(r+1)^n$) and Theorem3.6(1), and we have
\begin{center}
$B_n^{(-2,\overbrace {\scriptstyle{0,\ldots,0}}^{r-1})}$=$\dfrac{2}{(r+1)!} \displaystyle \sum_{k=1}^{r+1}
\begin{bmatrix}
r+1\\k
\end{bmatrix}
B_n^{(-k)}-\dfrac{1}{r!} \displaystyle \sum_{k=1}^{r}
\begin{bmatrix}
r\\k
\end{bmatrix}
B_n^{(-k)}$.
\end{center}

We consider similarly in the case of $m=3$, that is, $B_n^{(-3,\overbrace {\scriptstyle{0,\ldots,0}}^{r-1})}$. 
Then we have
\begin{center}
$B_n^{(-3,\overbrace {\scriptstyle{0,\ldots,0}}^{r-1})}$=$\dfrac{6}{(r+2)!} \displaystyle \sum_{k=1}^{r+2}
\begin{bmatrix}
r+2\\k
\end{bmatrix}
B_n^{(-k)}-\dfrac{6}{(r+1)!} \displaystyle \sum_{k=1}^{r+1}
\begin{bmatrix}
r+1\\k
\end{bmatrix}
B_n^{(-k)}$

\hspace{-4cm} $+\dfrac{1}{r!} \displaystyle \sum_{k=1}^{r}
\begin{bmatrix}
r\\k
\end{bmatrix}
B_n^{(-k)}$.
\end{center}
From here we can be considered the generalizations, that is, $B_n^{(-m,\overbrace {\scriptstyle{0,\ldots,0}}^{r-1})}$ as follows.

\textbf{Theorem3.9.} We have the following relations on  $B_n^{(-m,\overbrace {\scriptstyle{0,\ldots,0}}^{r-1})}$
\begin{center}
$B_n^{(-m,\overbrace {\scriptstyle{0,\ldots,0}}^{r-1})}$=$\displaystyle \sum_{l=1}^{m}
\dfrac{(-1)^{l+m} l! \genfrac{\{}{\}}{0pt}{}{m}{l}}{(r+l-1)!} \displaystyle \sum_{k=1}^{r+l-1}
\begin{bmatrix}
r+l-1\\k
\end{bmatrix}
B_n^{(-k)}$.
\end{center}

Proof. Using Theorem2.6(1) and Theorem3.6(1), we have
\begin{center}
$ \dfrac{1}{r!} \displaystyle \sum_{k=1}^{r}
\begin{bmatrix}
r\\k
\end{bmatrix}
B_n^{(-k)}$=$(r+1)^n$.
\end{center}

Here, we replace $r \to r+l-1$, and we have
\begin{center}
$ \dfrac{1}{(r+l-1)!} \displaystyle \sum_{k=1}^{r+l-1}
\begin{bmatrix}
r+l-1\\k
\end{bmatrix}
B_n^{(-k)}$=$(r+l)^n$.
\end{center}

Hence we obtain
\begin{center}
$\displaystyle \sum_{l=1}^{m}
\dfrac{(-1)^{l+m} l! \genfrac{\{}{\}}{0pt}{}{m}{l}}{(r+l-1)!} \displaystyle \sum_{k=1}^{r+l-1}
\begin{bmatrix}
r+l-1\\k
\end{bmatrix}
B_n^{(-k)}$
=$\displaystyle \sum_{l=1}^{m} (-1)^{l+m} l! \genfrac{\{}{\}}{0pt}{}{m}{l} (r+l)^n$.
\end{center}
Here the right hand of the last equality can be obtained by putting $k_1=m$,\hspace{3pt} $k_2=\cdots=k_r=0$ in Theorem2.1.
Thus it equals to $B_n^{(-m,\overbrace {\scriptstyle{0,\ldots,0}}^{r-1})}$ and we obtain Theorem3.9. 
\hspace{12.5cm} $\Box$

\vspace{5pt}
By using Lemma3.4, we can also express Theorem3.9 by using the Stirling numbers of the second kind.

Next, we consider $B_n^{\overbrace {\scriptstyle{(0,\ldots,0}}^{r-1},-m)}$ which are the generalizations of Theorem3.6(2) in the same way. We consider the small values on $m$, and we have

$B_n^{\overbrace {\scriptstyle{(0,\ldots,0}}^{r-1},-1)}$
=$\genfrac{\{}{\}}{0pt}{}{1}{1}r(r+1)^n$ 

$B_n^{\overbrace {\scriptstyle{(0,\ldots,0}}^{r-1},-2)}$
=$-\genfrac{\{}{\}}{0pt}{}{2}{1}r(r+1)^n+\genfrac{\{}{\}}{0pt}{}{2}{2}r(r+1)(r+2)^n$

$B_n^{\overbrace {\scriptstyle{(0,\ldots,0}}^{r-1},-3)}$
=$\genfrac{\{}{\}}{0pt}{}{3}{1}r(r+1)^n-\genfrac{\{}{\}}{0pt}{}{3}{2}r(r+1)(r+2)^n
+\genfrac{\{}{\}}{0pt}{}{3}{3}r(r+1)(r+2)(r+3)^n$

\begin{center}
$\cdots$.
\end{center}
From here we can be considered the generalizations, that is, $B_n^{\overbrace {\scriptstyle{(0,\ldots,0}}^{r-1},-m)}$ as follows.

\textbf{Theorem3.10.} We have the folowing relations on $B_n^{\overbrace {\scriptstyle{(0,\ldots,0}}^{r-1},-m)}$
\begin{center}
$B_n^{\overbrace {\scriptstyle{(0,\ldots,0}}^{r-1},-m)}$
=$\displaystyle \sum_{l=1}^{m} \dfrac{(-1)^{l+m} (r)_{l-1} \genfrac{\{}{\}}{0pt}{}{m}{l} }
{(r+l-2)!} \displaystyle \sum_{k=1}^{r+l-1}
\begin{bmatrix}
r+l-1\\k
\end{bmatrix}
B_n^{(-k)}$.
\end{center}
Here, we define $(r)_{l}=r(r+1)\cdots (r+l-1)$ and $(r)_0=1$.

\vspace{5pt}
Proof. Using Theorem2.6(2) and Theorem3.6(2), we have
\begin{center}
$\dfrac{1}{(r-1)!} \displaystyle \sum_{k=1}^{r}
\begin{bmatrix}
r\\k
\end{bmatrix}
B_n^{(-k)}$=$r(r+1)^n$.
\end{center}
Here we replace $r \to r+l-1$, and we have
\begin{center}
$\dfrac{1}{(r+l-2)!} \displaystyle \sum_{k=1}^{r+l-1}
\begin{bmatrix}
r+l-1\\k
\end{bmatrix}
B_n^{(-k)}$=$(r+l-1)(r+l)^n$.
\end{center}
Hence we obtain
\begin{align*}
\displaystyle \sum_{l=1}^{m}\dfrac{(-1)^{l+m}(r)_{l-1}\genfrac{\{}{\}}{0pt}{}{m}{l}}{(r+l-2)!}
\displaystyle \sum_{k=1}^{r+l-1}
\begin{bmatrix}
r+l-1\\k
\end{bmatrix}
B_n^{(-k)}
&=\displaystyle \sum_{l=1}^{m} (-1)^{l+m} (r)_{l-1} \genfrac{\{}{\}}{0pt}{}{m}{l} (r+l-1)(r+l)^n\\
&=\displaystyle \sum_{l=1}^{m} (-1)^{l+m} (r)_{l} \genfrac{\{}{\}}{0pt}{}{m}{l} (r+l)^n.
\end{align*}
Here from Theorem2.1, since we have $B_n^{\overbrace {\scriptstyle{(0,\ldots,0}}^{r-1},-m)}$=$\displaystyle \sum_{l=1}^{m} {\alpha}_l^{(0,\ldots,0,m)} (l+r)^n$,
we prove by induction on $m$

\vspace{-5pt}
and $l$ that we have $\alpha_l^{\overbrace {\scriptstyle{(0,\ldots,0}}^{r-1},m)}$=$(-1)^{l+m}(r)_l  \genfrac{\{}{\}}{0pt}{}{m}{l}$.

\vspace{5pt}
\hspace{-15pt} First, we prove that we have $\alpha_1^{(0,\ldots,0,m)}$=$(-1)^{m+1} r \cdots (B)$.

\hspace{-15pt} From Theorem2.1(iii), since we have $\alpha_1^{(0,\ldots,0,m)}=r\alpha_0^{(0,\ldots,0,m-1)}-\alpha_1^{(0,\ldots,0,m-1)}$,
we have (B) for $m=1$.

\hspace{-15pt} We assume that we have $\alpha_1^{(0,\ldots,0,k)}$=$(-1)^{k+1}r$ for $m=k$ ($k\geq 1$).

\hspace{-15pt} For $m=k+1$, we have
\begin{align*}
\alpha_1^{(0,\ldots,0,k+1)}&=r\alpha_0^{(0,\ldots,0,k)}-\alpha_1^{(0,\ldots,0,k)}\\
&=-\alpha_1^{(0,\ldots,0,k)}\\
&=-(-1)^{k+1} r\\
&=(-1)^{k+2} r.
\end{align*}
Since this shows that (B) is true for $m=k+1$, we have (B) for all integers $m$.

\hspace{-15pt} Next,we prove that we have $\alpha_l^{(0,\ldots,0,m)}$=$(-1)^{l+m} (r)_l \genfrac{\{}{\}}{0pt}{}{m}{l} \cdots (C)$.

\hspace{-15pt} From (B), we have (C) for $m=1$.

\hspace{-15pt} We assume that we have $\alpha_k^{(0,\ldots,0,m)}$=$(-1)^{k+m} (r)_k\genfrac{\{}{\}}{0pt}{}{m}{k}$ for $l=k$ ($k \geq 1$).

\hspace{-15pt} For $l=k+1$, we have
\begin{align*}
\alpha_{k+1}^{(0,\ldots,0,m)}&=(k+r)\alpha_k^{(0,\ldots,0,m-1)}-(k+1)\alpha_{k+1}^{(0,\ldots,0,m-1)}\\
&=(k+r)\alpha_k^{(0,\ldots,0,m-1)}-(k+1)  {\{(k+r)\alpha_k^{(0,\ldots,0,m-2)}-(k+1)\alpha_{k+1}^{(0,\ldots,0,m-2)} }\}\\
&=(k+r)\alpha_k^{(0,\ldots,0,m-1)}-(k+1)(k+r) \alpha_k^{(0,\ldots,0,m-2)}+(k+1)^2  \alpha_{k+1}^{(0,\ldots,0,m-2)}\\
&=(-1)^{k+m-1} (r)_k (k+r) \genfrac{\{}{\}}{0pt}{}{m-1}{k}-(-1)^{k+m-2} (r)_k (k+1)(k+r) \genfrac{\{}{\}}{0pt}{}{m-2}{k}\\
& \hspace{0.4cm} +(k+1)^2  \alpha_{k+1}^{(0,\ldots,0,m-2)}\\
&\ldots\\
&=(-1)^{k+m-1} (r)_{k+1} \biggl[\genfrac{\{}{\}}{0pt}{}{m-1}{k}+(k+1)\genfrac{\{}{\}}{0pt}{}{m-2}{k}+ \cdots +
(k+1)^{m-k-1}\\
& \hspace{0.45cm} \times \genfrac{\{}{\}}{0pt}{}{m-(m-k)}{k} \biggr].
\end{align*}
Since the Stirling numbers of the second kind satisfy the recurrence formula
\begin{center}
$\genfrac{\{}{\}}{0pt}{}{n+1}{m+1}=\genfrac{\{}{\}}{0pt}{}{n}{m}+(m+1)\genfrac{\{}{\}}{0pt}{}{n}{m+1}$,
\end{center}

\hspace{-10pt} We have  $\genfrac{\{}{\}}{0pt}{}{m-1}{k}+(k+1)\genfrac{\{}{\}}{0pt}{}{m-2}{k}+(k+1)^2\genfrac{\{}{\}}{0pt}{}{m-3}{k}+ \cdots +
(k+1)^{m-k-1} \genfrac{\{}{\}}{0pt}{}{m-(m-k)}{k}$

$\hspace{1cm}=\displaystyle \sum_{i=1}^{m-k} (k+1)^{i-1} \genfrac{\{}{\}}{0pt}{}{m-i}{k}$

$\hspace{1cm}=\displaystyle \sum_{i=1}^{m-k} (k+1)^{i-1} \bigg( \genfrac{\{}{\}}{0pt}{}{m-i+1}{k+1} - (k+1)\genfrac{\{}{\}}{0pt}{}{m-i}{k+1}\biggr)$

$\hspace{1cm}=\displaystyle \sum_{i=1}^{m-k} \bigg[(k+1)^{i-1} \genfrac{\{}{\}}{0pt}{}{m-i+1}{k+1} - (k+1)^i\genfrac{\{}{\}}{0pt}{}{m-i}{k+1}\biggr]$

$\hspace{1cm}=\genfrac{\{}{\}}{0pt}{}{m}{k+1}-(k+1) \genfrac{\{}{\}}{0pt}{}{m-1}{k+1}+(k+1) \genfrac{\{}{\}}{0pt}{}{m-1}{k+1}- (k+1)^2 \genfrac{\{}{\}}{0pt}{}{m-2}{k+1}+ \cdots$

\hspace{1.2cm} $+(k+1)^{m-k-1}  \genfrac{\{}{\}}{0pt}{}{k+1}{k+1}-(k+1)^{m-k} \genfrac{\{}{\}}{0pt}{}{k}{k+1}$

$\hspace{1cm}=\genfrac{\{}{\}}{0pt}{}{m}{k+1}$.

\hspace{-15pt} Therefore we have $\alpha_{k+1}^{(0,\ldots,0,m)}$=$(-1)^{k+m-1} (r)_{k+1} \genfrac{\{}{\}}{0pt}{}{m}{k+1}$ and  this shows that (C) is true for $l=k+1$.

\hspace{-15pt} Hence, we have (C) for all integers $l,m$, and this completes the proof. \hspace{1.8cm} ${\Box}$

By using Lemma3.4, we can also express Theorem3.10 by using the Stirling numbers of the second kind.

Finally, we consider the case of $B_n^{\overbrace{\scriptstyle{(0,\ldots,0,-m,0,\ldots,0)}}^r}$ (where i-th component is $-m$
and others are $0$) which are the extension of
Theorem3.9 and Theorem3.10. 
For example, we fluctuate the value of 2-th. Then we obtain 

\hspace{4cm} $B_n^{\overbrace{\scriptstyle{(0,-1,0,\ldots,0)}}^r}$=$2(r+1)^n$

\hspace{4cm} $B_n^{\overbrace{\scriptstyle{(0,-2,0,\ldots,0)}}^r}$=$-2(r+1)^n+6(r+2)^n$

\hspace{4cm} $B_n^{\overbrace{\scriptstyle{(0,-3,0,\ldots,0)}}^r}$=$2(r+1)^n-6(r+2)^n+12(r+3)^n$

\hspace{4cm} $B_n^{\overbrace{\scriptstyle{(0,-4,0,\ldots,0)}}^r}$=$-2(r+1)^n+6(r+2)^n-12(r+3)^n+24(r+4)^n$

\begin{center}
$\cdots$.
\end{center}
By using this relations and the recurrence formula on Theorem2.1, we can be considered the following relations.

\textbf{Theorem3.11.} We have the following relations on $B_n^{\overbrace{\scriptstyle{(0,\ldots,0,-m,0,\ldots,0)}}^r}$,
and i-th component is $-m$ and others are $0$;
\begin{center}
$B_n^{\overbrace{\scriptstyle{(0,\ldots,0,-m,0,\ldots,0)}}^r}$=$\displaystyle \sum_{l=1}^{m}
\dfrac{(-1)^{l-m} (i)_l \genfrac{\{}{\}}{0pt}{}{m}{l} }{(r+l-1)!} 
\displaystyle \sum_{k=1}^{r+l-1}
\begin{bmatrix}
r+l-1\\k
\end{bmatrix}B_n^{(-k)}$.
\end{center}

Proof. In the proof of Theorem3.9, we have
\begin{center}
$\dfrac{1}{(r+l-1)!} 
\displaystyle \sum_{k=1}^{r+l-1}
\begin{bmatrix}
r+l-1\\k
\end{bmatrix}B_n^{(-k)}$
=$(r+l)^n$.
\end{center}
Futhermore, using Theorem2.1 and $\alpha_l^{\overbrace {\scriptstyle{(0,\ldots,0}}^{r-1},m)}$=$(-1)^{l+m} (r)_l\genfrac{\{}{\}}{0pt}{}{m}{l}$ in the proof of Theorem3.10, we have
\begin{align*}
\displaystyle \sum_{l=1}^{m}\dfrac{(-1)^{l-m} (i)_l \genfrac{\{}{\}}{0pt}{}{m}{l} }{(r+l-1)!} 
\displaystyle \sum_{k=1}^{r+l-1}
\begin{bmatrix}
r+l-1\\k
\end{bmatrix}B_n^{(-k)}
&=\displaystyle \sum_{l=1}^{m} (-1)^{l-m} (i)_l \genfrac{\{}{\}}{0pt}{}{m}{l} (r+l)^n\\
&=\displaystyle \sum_{l=1}^{m} \alpha_l^{\overbrace {\scriptstyle{(0,\ldots,0}}^{i-1},m)} (r+l)^n\\
&=\displaystyle \sum_{l=1}^{m} \alpha_l^{\overbrace {\scriptstyle{(0,\ldots,0}}^{i-1},m,0,\ldots,0)}(r+l)^n\\
&=B_n^{\overbrace{\scriptstyle{(0,\ldots,0,-m,0,\ldots,0)}}^r}.
\end{align*}
Hence we obtain
\begin{center}
$B_n^{\overbrace{\scriptstyle{(0,\ldots,0,-m,0,\ldots,0)}}^r}$=$\displaystyle \sum_{l=1}^{m}
\dfrac{(-1)^{l+m} (i)_l \genfrac{\{}{\}}{0pt}{}{m}{l} }{(r+l-1)!} 
\displaystyle \sum_{k=1}^{r+l-1}
\begin{bmatrix}
r+l-1\\k
\end{bmatrix}{B_n^{(-k)}}$,
\end{center}
and This completes the proof.  \hspace{8.5cm} ${\Box}$

By using Lemma3.4, we can also express Theorem3.11 by using the Stirling numbers of the second kind.

If we put $i=1$, $i=r$ in Theorem3.11, we obtain the following
\begin{center}
$B_n^{(-m,\overbrace {\scriptstyle{0,\ldots,0}}^{r-1})}$=$\displaystyle \sum_{l=1}^{m}
\dfrac{(-1)^{l+m} l! \genfrac{\{}{\}}{0pt}{}{m}{l}}{(r+l-1)!} \displaystyle \sum_{k=1}^{r+l-1}
\begin{bmatrix}
r+l-1\\k
\end{bmatrix}
{B_n^{(-k)}}$,
\end{center}

\begin{center}
\hspace{0.5cm} $B_n^{\overbrace {\scriptstyle{(0,\ldots,0}}^{r-1},-m)}$
=$\displaystyle \sum_{l=1}^{m} \dfrac{(-1)^{l+m} (r)_{l-1} \genfrac{\{}{\}}{0pt}{}{m}{l} }
{(r+l-2)!} \displaystyle \sum_{k=1}^{r+l-1}
\begin{bmatrix}
r+l-1\\k
\end{bmatrix}
{B_n^{(-k)}}_{\textbf.}$
\end{center}
Hence we find that Theorem3.11 is the extension of Theorem3.9 and Theorem3.10.

Here we represent $B_n^{(-m,\overbrace {\scriptstyle{0,\ldots,0}}^{r-1})}$,\,$B_n^{\overbrace {\scriptstyle{(0,\ldots,0}}^{r-1},-m)}$,\,$B_n^{\overbrace{\scriptstyle{(0,\ldots,0,-m,0,\ldots,0)}}^r}$ in the form of powers on  $r+l$ ($1 \leq  l \leq m $), and we see the sum of coefficients.
For example, we put $m=1$. Then we have the following relations from Theorem2.6:

\hspace{4cm} $B_n^{(-1,\overbrace {\scriptstyle{0,\ldots,0}}^{r-1})}=(r+1)^n$,

\hspace{4.05cm} $B_n^{\overbrace {\scriptstyle{(0,\ldots,0}}^{r-1},-1)}=r(r+1)^n$,

\hspace{3.25cm} $B_n^{\overbrace{\scriptstyle{(0,\ldots,0,-1,0,\ldots,0)}}^r}=i(r+1)^n$.

Hence each coefficients are $1$,\,$r$, and $i$. From this results, we can be considered the following.

\textbf{Theorem3.12.} We have the following relations on the sum of coefficients

(1)\hspace{3pt} The sum of coefficients on $B_n^{(-m,\overbrace {\scriptstyle{0,\ldots,0}}^{r-1})}$ are $1$.

(2)\hspace{3pt} The sum of coefficients on $B_n^{\overbrace {\scriptstyle{(0,\ldots,0}}^{r-1},-m)}$ are $r^m$.

(3)\hspace{3pt} The sum of coefficients on $B_n^{\overbrace{\scriptstyle{(0,\ldots,0,-m,0,\ldots,0)}}^r}$ are $i^m$.

Proof. (1) In the proof of Theorem3.9, we have
\begin{center}
$B_n^{(-m,\overbrace {\scriptstyle{0,\ldots,0}}^{r-1})}$
=$\displaystyle \sum_{l=1}^{m} (-1)^{l+m} l! \genfrac{\{}{\}}{0pt}{}{m}{l} (r+l)^n$.
\end{center}
Hence it suffices to show that
$\displaystyle \sum_{l=1}^{m}(-1)^{l+m} l! \genfrac{\{}{\}}{0pt}{}{m}{l}=1$. We have
\begin{align*}
\hspace{2.5cm} \displaystyle \sum_{l=1}^{m}(-1)^{l+m} l! \genfrac{\{}{\}}{0pt}{}{m}{l}
&=(-1)^m \displaystyle \sum_{l=1}^{m} (-1)^l \displaystyle \sum_{k=0}^{l}
\begin{bmatrix}
l\\k
\end{bmatrix}
\genfrac{\{}{\}}{0pt}{}{m}{l}\\
&=(-1)^m \displaystyle \sum_{l=0}^{m} (-1)^l \displaystyle \sum_{k=0}^{m}
\begin{bmatrix}
l\\k
\end{bmatrix}
\genfrac{\{}{\}}{0pt}{}{m}{l}\\
&=(-1)^m \displaystyle \sum_{k=0}^{m}\displaystyle \sum_{l=0}^{m} (-1)^l 
\genfrac{\{}{\}}{0pt}{}{m}{l}
\begin{bmatrix}
l\\k
\end{bmatrix}.
\end{align*}
Here since $\displaystyle \sum_{l=0}^{n} (-1)^l \genfrac{\{}{\}}{0pt}{}{n}{l}
\begin{bmatrix}
l\\m
\end{bmatrix}
=(-1)^m \delta_{m,n}$ ([1]), we obtain
\begin{align*}
\displaystyle \sum_{l=1}^{m}(-1)^{l+m} l! \genfrac{\{}{\}}{0pt}{}{m}{l}
&=(-1)^m \displaystyle \sum_{k=0}^{m} (-1)^k \delta_{k,m}\\
&=(-1)^m (-1)^m \delta_{m,m}\\
&=1.
\end{align*}
(2) In the proof of Theorem3.10, we have
\begin{center}
$B_n^{\overbrace {\scriptstyle{(0,\ldots,0}}^{r-1},-m)}$=$\displaystyle \sum_{l=1}^{m} (-1)^{l+m} (r)_l \genfrac{\{}{\}}{0pt}{}{m}{l} (r+l)^n$.
\end{center}
Hence it suffices to show that 
$\displaystyle \sum_{l=1}^{m} (-1)^{l+m} (r)_l \genfrac{\{}{\}}{0pt}{}{m}{l}=r^m$.

Since $x^n=\displaystyle \sum_{k=1}^{n} \genfrac{\{}{\}}{0pt}{}{n}{k} (-1)^{n-k} (x)_k$\hspace{1pt}
($n \geq 0$) ([6]), 
$\displaystyle \sum_{l=1}^{m} (-1)^{l+m} (r)_l \genfrac{\{}{\}}{0pt}{}{m}{l}=r^m$.

\vspace{5pt}
(3) In the proof of Theorem3.11, we have
\begin{center}
$B_n^{\overbrace{\scriptstyle{(0,\ldots,0,-m,0,\ldots,0)}}^r}$=$\displaystyle \sum_{l=1}^{m} (-1)^{l+m} (i)_l \genfrac{\{}{\}}{0pt}{}{m}{l} (r+l)^n$.
\end{center}
Hence it suffices to show that 
$\displaystyle \sum_{l=1}^{m} (-1)^{l+m} (i)_l \genfrac{\{}{\}}{0pt}{}{m}{l}=i^m$.
This identity can be obtained by putting $r=i$ in (2), and we obtain the result.
\hspace{4cm} ${\Box}$

For $m=1$, since each coefficients are $1$,\,$r$, and $i$, we find that Theorem3.12 is the generalizations.

We introduced several relations between Multi-Poly-Bernoulli numbers and Poly-Bernoulli numbers up to here.
We obtained Theorem3.1 by fluctuating values of $r$ which represent numbers of $0$. Here we consider that fluctuating values of $m$ which represent numbers except for $0$.

\textbf{Theorem3.13.} We have the following relations

(1)\hspace{3pt} $B_n^{(-k)}=\displaystyle \sum_{m=0}^{k-1} (-1)^{k-m-1} \binom{k}{m} B_n^{(-m,0)} \hspace{5pt} (k \geq 1)$.

\hspace{10pt} Where the sum of coefficients on $B_n^{(-m,0)}$ are $1$.

(2)\hspace{3pt} $B_n^{(-k)}=\displaystyle \sum_{m=0}^{k-1} (-1)^{k-m-1} \binom{k-1}{m} B_n^{(0,-m)} \hspace{5pt} (k \geq 1)$.

\hspace{10pt} Where the sum of coefficients on $B_n^{(0,-m)}$ are $0$ for $k \geq 2$. we regard the sum on the right hand as $1$ if $k=1$.

\begin{align*}
\mathrm{Proof.} (1)\, \mathrm{R.H.S.} &=\displaystyle \sum_{m=0}^{k-1} (-1)^{k-m-1} \binom{k}{m}
\displaystyle \sum_{l=1}^{m} (-1)^{l+m} l! \genfrac{\{}{\}}{0pt}{}{m}{l} (l+2)^n\\
&=\displaystyle \sum_{m=0}^{k-1} \binom{k}{m} \displaystyle \sum_{l=1}^{m} (-1)^{l+1+k}
l! \genfrac{\{}{\}}{0pt}{}{m}{l} (l+2)^n\\
&=\displaystyle \sum_{m=0}^{k-1} \binom{k}{m} \displaystyle \sum_{l=0}^{m} (-1)^{l+1+k}
l! \genfrac{\{}{\}}{0pt}{}{m}{l} (l+2)^n\\
&=\displaystyle \sum_{m=0}^{k-1} \binom{k}{m} \displaystyle \sum_{l=0}^{k-1} (-1)^{l+1+k}
l! \genfrac{\{}{\}}{0pt}{}{m}{l} (l+2)^n\\
&=\displaystyle \sum_{l=0}^{k-1} (-1)^{l+1+k} l! (l+2)^n \displaystyle \sum_{m=0}^{k-1}
\binom{k}{m} \genfrac{\{}{\}}{0pt}{}{m}{l}\\
&=\displaystyle \sum_{l=0}^{k-1} (-1)^{l+1+k} l! (l+2)^n \biggl[\displaystyle \sum_{m=0}^{k}
\binom{k}{m} \genfrac{\{}{\}}{0pt}{}{m}{l}-\genfrac{\{}{\}}{0pt}{}{k}{l}\biggr]\\
&=\displaystyle \sum_{l=0}^{k-1} (-1)^{l+1+k} l! (l+2)^n \biggl[\genfrac{\{}{\}}{0pt}{}{k+1}{l+1}-\genfrac{\{}{\}}{0pt}{}{k}{l}\biggr]\\
&=\displaystyle \sum_{l=0}^{k-1} (-1)^{l+1+k} l! (l+2)^n (l+1)  \genfrac{\{}{\}}{0pt}{}{k}{l+1}\\
&=\displaystyle \sum_{l=0}^{k-1} (-1)^{l+1+k} (l+1)! \genfrac{\{}{\}}{0pt}{}{k}{l+1} (l+2)^n\\
&=\displaystyle \sum_{l=1}^{k} (-1)^{l+k} l! \genfrac{\{}{\}}{0pt}{}{k}{l} (l+1)^n\\
&=B_n^{(-k)}.
\end{align*}
Hence we obtain the result. Futhermore, the sum of coefficients on $B_n^{(-m,0)}$ are
\begin{align*}
\displaystyle \sum_{m=0}^{k-1} (-1)^{k-m-1} \binom{k}{m}
&=\displaystyle \sum_{m=0}^{k} (-1)^{k-m-1} \binom{k}{m}+1\\
\end{align*}
\begin{align*}
\hspace{3cm}&=-\displaystyle \sum_{m=0}^{k} (-1)^{k-m} \binom{k}{m}+1\\
&=-(1-1)^k+1\\
&=1.
\end{align*}

This completes the proof and we obtain (1).

(2)We put $r=2$ in Theorem3.10, and we use Theorem3.6 and Theorem2.6. Then we have
\begin{align*} B_n^{(0,-m)}&=\displaystyle \sum_{l=1}^{m} \dfrac{(-1)^{l+m} (2)_{l-1} 
\genfrac{\{}{\}}{0pt}{}{m}{l}}{l!} \displaystyle \sum_{k=1}^{l+1}
\begin{bmatrix}
l+1\\k
\end{bmatrix}B_n^{(-k)}\\
&=\displaystyle \sum_{l=1}^{m} (-1)^{l+m} l! \genfrac{\{}{\}}{0pt}{}{m}{l}
\dfrac{1}{l!} \displaystyle \sum_{k=1}^{l+1}
\begin{bmatrix}
l+1\\k
\end{bmatrix}B_n^{(-k)}\\
&=\displaystyle \sum_{l=1}^{m} (-1)^{l+m} l! \genfrac{\{}{\}}{0pt}{}{m}{l}
B_n^{\overbrace {\scriptstyle{(0,\cdots,0}}^{l},-1)}\\
&=\displaystyle \sum_{l=1}^{m} (-1)^{l+m} l! \genfrac{\{}{\}}{0pt}{}{m}{l} 
(l+1) (l+2)^n. 
\end{align*}
Hence we substitute this identity on the right of Theorem3.13, and we have
\begin{align*}
\mathrm{R.H.S.} &=\displaystyle \sum_{m=0}^{k-1} (-1)^{k-m-1} \binom{k-1}{m}
\displaystyle \sum_{l=1}^{m} (-1)^{l+m} l! \genfrac{\{}{\}}{0pt}{}{m}{l} 
(l+1) (l+2)^n\\
&=\displaystyle \sum_{m=0}^{k-1} \binom{k-1}{m} \displaystyle \sum_{l=1}^{m}
(-1)^{l+1+k} (l+1)! \genfrac{\{}{\}}{0pt}{}{m}{l} 
(l+2)^n\\
&=\displaystyle \sum_{m=0}^{k-1} \binom{k-1}{m} \displaystyle \sum_{l=0}^{m}
(-1)^{l+1+k} (l+1)! \genfrac{\{}{\}}{0pt}{}{m}{l} 
(l+2)^n\\
&=\displaystyle \sum_{m=0}^{k-1} \binom{k-1}{m} \displaystyle \sum_{l=0}^{k-1}
(-1)^{l+1+k} (l+1)! \genfrac{\{}{\}}{0pt}{}{m}{l} 
(l+2)^n\\
&=\displaystyle \sum_{l=0}^{k-1} (-1)^{l+1+k} (l+1)! (l+2)^n
\displaystyle \sum_{m=0}^{k-1} \binom{k-1}{m} \genfrac{\{}{\}}{0pt}{}{m}{l}\\
&=\displaystyle \sum_{l=0}^{k-1} (-1)^{l+1+k} (l+1)! (l+2)^n \genfrac{\{}{\}}{0pt}{}{k}{l+1}\\
\end{align*}
\begin{align*}
\hspace{-2cm}&=\displaystyle \sum_{l=0}^{k-1} (-1)^{l+1+k} (l+1)! \genfrac{\{}{\}}{0pt}{}{k}{l+1} (l+2)^n\\
&=\displaystyle \sum_{l=1}^{k} (-1)^{l+k} l! \genfrac{\{}{\}}{0pt}{}{k}{l} (l+1)^n\\
&=B_n^{(-k)}.
\end{align*}
Therefore we obtain the result. Futhermore, the sum of coefficients on $B_n^{(0,-m)}$ are
\begin{align*}
\displaystyle \sum_{m=0}^{k-1} (-1)^{k-m-1} \binom{k-1}{m}
&=(1-1)^{k-1}=\begin{cases}\
1 & \text{($k=1$)}\\
\hspace{3pt} 0 & \text{($k \geq 2$)}.
\end{cases}
\end{align*}
This completes the proof and we obtain (2). \hspace{6cm} ${\Box}$

We extend Theorem3.13, and we can write $B_n^{(-k)}$ by using the sum of $B_n^{\overbrace {\scriptstyle{(0,\ldots,0}}^{l+1})}$ \hspace{3pt}
($ l \geq 1$) and them of $B_n^{(-m,\overbrace {\scriptstyle{0,\ldots,0}}^{r})}$, or
$B_n^{\overbrace {\scriptstyle{(0,\ldots,0}}^{r},-m)}$. First, we write $B_n^{(-k)}$ by using the sum of $B_n^{\overbrace {\scriptstyle{(0,\ldots,0}}^{l+1})}$ \hspace{3pt}
($ l \geq 1$) and them of $B_n^{(-m,\overbrace {\scriptstyle{0,\ldots,0}}^{r})}$.

\textbf{Theorem3.14.} We have the following relations

$B_n^{(-k)}=\displaystyle \sum_{l=1}^{r} (-1)^{k-l} l! \genfrac{\{}{\}}{0pt}{}{k}{l}
B_n^{\overbrace {\scriptstyle{(0,\ldots,0}}^{l+1})}+\displaystyle \sum_{m=1}^{k-r}
(-1)^{k-m-r}  \binom{k}{m} r! \genfrac{\{}{\}}{0pt}{}{k-m}{r} B_n^{(-m,\overbrace{\scriptstyle {0,\ldots,0}}^{r})}$.

Proof. (1)\,$\displaystyle \sum_{l=1}^{r} (-1)^{k-l} l! \genfrac{\{}{\}}{0pt}{}{k}{l}
B_n^{\overbrace {\scriptstyle{(0,\ldots,0}}^{l+1})}=\displaystyle \sum_{l=1}^{r}
(-1)^{k-l} l! \genfrac{\{}{\}}{0pt}{}{k}{l} (l+1)^n$.

Futhermore, we have

\hspace{7pt} $\displaystyle \sum_{m=1}^{k-r} (-1)^{k-m-r}  \binom{k}{m} r! \genfrac{\{}{\}}{0pt}{}{k-m}{r} B_n^{(-m,\overbrace {\scriptstyle{0,\ldots,0}}^{r})}$

$=\displaystyle \sum_{m=1}^{k-r} (-1)^{k-m-r}  \binom{k}{m} r! \genfrac{\{}{\}}{0pt}{}{k-m}{r} \displaystyle \sum_{l=1}^{m}
(-1)^{l+m} l! \genfrac{\{}{\}}{0pt}{}{m}{l} (l+r+1)^n$

$=\displaystyle \sum_{m=r+1}^{k} (-1)^{k-m}  \binom{k}{m-r} r! \genfrac{\{}{\}}{0pt}{}{k-m+r}{r} \displaystyle \sum_{l=1}^{m-r}
(-1)^{l+m-r} l! \genfrac{\{}{\}}{0pt}{}{m-r}{l} (l+r+1)^n$

$=\displaystyle \sum_{m=r+1}^{k} \binom{k}{m-r} r! \genfrac{\{}{\}}{0pt}{}{k-m+r}{r} \displaystyle \sum_{l=1}^{k-1}
(-1)^{k-l-r} l! \genfrac{\{}{\}}{0pt}{}{m-r}{l} (l+r+1)^n$

$=\displaystyle \sum_{m=r+1}^{k} \binom{k}{m-r} r! \genfrac{\{}{\}}{0pt}{}{k-m+r}{r} \displaystyle \sum_{l=1}^{k}
(-1)^{k-l-r} l! \genfrac{\{}{\}}{0pt}{}{m-r}{l} (l+r+1)^n$

$=\displaystyle \sum_{l=1}^{k} (-1)^{k-l-r} r! l! (l+r+1)^n \displaystyle \sum_{m=r+1}^{k} \genfrac{\{}{\}}{0pt}{}{m-r}{l} \genfrac{\{}{\}}{0pt}{}{k-m+r}{r} \binom{k}{m-r}$

$=\displaystyle \sum_{l=1}^{k} (-1)^{k-l-r} r! l! (l+r+1)^n \genfrac{\{}{\}}{0pt}{}{k}{l+r} \binom{l+r}{l}$

$=\displaystyle \sum_{l=1}^{k-r} (-1)^{k-l-r} r! l! (l+r+1)^n \genfrac{\{}{\}}{0pt}{}{k}{l+r} \binom{l+r}{l}$

$=\displaystyle \sum_{l=1}^{k-r} (-1)^{k-l-r} \genfrac{\{}{\}}{0pt}{}{k}{l+r} (l+r+1)^n r! l! \dfrac{(l+r)!}{l! r!}$

$=\displaystyle \sum_{l=1}^{k-r} (-1)^{k-l-r} (l+r)! \genfrac{\{}{\}}{0pt}{}{k}{l+r}
(l+r+1)^n$

$=\displaystyle \sum_{l=r+1}^{k} (-1)^{k-l} l! \genfrac{\{}{\}}{0pt}{}{k}{l}
(l+1)^n$ \hspace{5pt} ($l \to l-r$).

Hence the right hand of Theorem3.14 can be expressed as follows, and we obtain the result.
\begin{align*}
\mathrm{R.H.S.} &=\displaystyle \sum_{l=1}^{r}
(-1)^{k-l} l! \genfrac{\{}{\}}{0pt}{}{k}{l} (l+1)^n+
\displaystyle \sum_{l=r+1}^{k} (-1)^{k-l} l! \genfrac{\{}{\}}{0pt}{}{k}{l}
(l+1)^n\\
&=\displaystyle \sum_{l=1}^{k} (-1)^{k-l} l! \genfrac{\{}{\}}{0pt}{}{k}{l} (l+1)^n\\
&=B_n^{(-k)}. \hspace{11cm} {\Box}
\end{align*}
\textbf{Example3.15.}

We give examples of Theorem3.14 for $1\leq r\leq2$,\,$1\leq k \leq 6$.

(i) For $r=1$

$B_n^{(-1)}=B_n^{(0,0)}$

$B_n^{(-2)}=-B_n^{(0,0)}+2B_n^{(-1,0)}$

$B_n^{(-3)}=B_n^{(0,0)}-3B_n^{(-1,0)}+3B_n^{(-2,0)}$

$B_n^{(-4)}=-B_n^{(0,0)}+4B_n^{(-1,0)}-6B_n^{(-2,0)}+4B_n^{(-3,0)}$

$B_n^{(-5)}=B_n^{(0,0)}-5B_n^{(-1,0)}+10B_n^{(-2,0)}-10B_n^{(-3,0)}+5B_n^{(-4,0)}$

$B_n^{(-6)}=-B_n^{(0,0)}+6B_n^{(-1,0)}-15B_n^{(-2,0)}+20B_n^{(-3,0)}-15B_n^{(-4,0)}+6B_n^{(-5,0)}$

(ii) For $r=2$

$B_n^{(-1)}=B_n^{(0,0)}$

$B_n^{(-2)}=-B_n^{(0,0)}+2B_n^{(0,0,0)}$

$B_n^{(-3)}=B_n^{(0,0)}-6B_n^{(0,0,0)}+6B_n^{(-1,0,0)}$

$B_n^{(-4)}=-B_n^{(0,0)}+14B_n^{(0,0,0)}-24B_n^{(-1,0,0)}+12B_n^{(-2,0,0)}$

$B_n^{(-5)}=B_n^{(0,0)}-30B_n^{(0,0,0)}+70B_n^{(-1,0,0)}-60B_n^{(-2,0,0)}+20B_n^{(-3,0,0)}$

$B_n^{(-6)}=-B_n^{(0,0)}+62B_n^{(0,0,0)}-180B_n^{(-1,0,0)}+210B_n^{(-2,0,0)}-120B_n^{(-3,0,0)}+30B_n^{(-4,0,0)}$

Next we write $B_n^{(-k)}$ by using the sum of $B_n^{\overbrace {\scriptstyle{(0,\ldots,0}}^{l+1})}$ \hspace{3pt}
($ l \geq 1$) and them of $B_n^{\overbrace {\scriptstyle{(0,\ldots,0}}^{r},-m)}$. First of all, we give concrete examples
for $1\leq r\leq3$,\,$1\leq k \leq 6$ such that Example3.15.

\textbf{Example3.16.}

We give examples for $1\leq r\leq3$,\,$1\leq k \leq 6$.

(i)For $r=1$

$B_n^{(-1)}=B_n^{(0,0)}$

$B_n^{(-2)}=-B_n^{(0,0)}+B_n^{(0,-1)}$

$B_n^{(-3)}=B_n^{(0,0)}-2B_n^{(0,-1)}+B_n^{(0,-2)}$

$B_n^{(-4)}=-B_n^{(0,0)}+3B_n^{(0,-1)}-3B_n^{(0,-2)}+B_n^{(0,-3)}$

$B_n^{(-5)}=B_n^{(0,0)}-4B_n^{(0,-1)}+6B_n^{(0,-2)}-4B_n^{(0,-3)}+B_n^{(0,-4)}$

$B_n^{(-6)}=-B_n^{(0,0)}+5B_n^{(0,-1)}-10B_n^{(0,-2)}+10B_n^{(0,-3)}-5B_n^{(0,-4)}+B_n^{(0,-5)}$

(ii)For $r=2$

$B_n^{(-1)}=B_n^{(0,0)}$

$B_n^{(-2)}=-B_n^{(0,0)}+2B_n^{(0,0,0)}$

$B_n^{(-3)}=B_n^{(0,0)}-6B_n^{(0,0,0)}+2B_n^{(0,0,-1)}$

$B_n^{(-4)}=-B_n^{(0,0)}+14B_n^{(0,0,0)}-10B_n^{(0,0,-1)}+2B_n^{(0,0,-2)}$

$B_n^{(-5)}=B_n^{(0,0)}-30B_n^{(0,0,0)}+34B_n^{(0,0,-1)}-14B_n^{(0,0,-2)}+2B_n^{(0,0,-3)}$

$B_n^{(-6)}=-B_n^{(0,0)} +62B_n^{(0,0,0)}-98B_n^{(0,0,-1)}+62B_n^{(0,0,-2)}-18B_n^{(0,0,-3)}+2B_n^{(0,0,-4)}$

(iii)For $r=3$

$B_n^{(-1)}=B_n^{(0,0)}$

$B_n^{(-2)}=-B_n^{(0,0)}+2B_n^{(0,0,0)}$

$B_n^{(-3)}=B_n^{(0,0)}-6B_n^{(0,0,0)}+6B_n^{(0,0,0,0)}$

$B_n^{(-4)}=-B_n^{(0,0)}+14B_n^{(0,0,0)}-36B_n^{(0,0,0,0)}+6B_n^{(0,0,0,-1)}$

$B_n^{(-5)}=B_n^{(0,0)}-30B_n^{(0,0,0)}+150B_n^{(0,0,0,0)}-54B_n^{(0,0,0,-1)}+6B_n^{(0,0,0,-2)}$

$B_n^{(-6)}=-B_n^{(0,0)}+62B_n^{(0,0,0)}-540B_n^{(0,0,0,0)}+312B_n^{(0,0,0,-1)}-72B_n^{(0,0,0,-2)}$

\hspace{1.4cm} $+6B_n^{(0,0,0,-3)}$

By Example3.16, we see the parts on coefficients except for plus or minus sign on Multi-Poly-Bernoulli\,numbers of the right hand. Then we obtain Pascal circles for $r=1$, and from here we can be considered the case of the generalizations as
follows.

\textbf{Conjecture3.17.} We will have the following relations

$B_n^{(-k)}=\displaystyle \sum_{l=1}^{r} (-1)^{k-l} l! \genfrac{\{}{\}}{0pt}{}{k}{l} 
B_n^{\overbrace {\scriptstyle{(0,\ldots,0}}^{l+1})}+
\displaystyle \sum_{m=1}^{k-r} (-1)^{k-m-r} a_{k-r-1,m} {B_n^{\overbrace {\scriptstyle{(0,\ldots,0}}^{r},-m)}}_{\textbf.}$

We define that $a_{k-r-1,m}=a_{k-r-2,m-1}+r a_{k-r-2,m}$ \hspace{3pt} ($0 \leq k-r-2,\,1 \leq m \leq k-r$),
 $a_{k-r-1,0}=r! \genfrac{\{}{\}}{0pt}{}{k}{r}$, and $a_{k-r-1,k-r}=r!$.

For example, for $k=4$ and $r=2$,
\vspace{-8pt}
\begin{align*}
B_n^{(-4)}&=\displaystyle \sum_{l=1}^{2} (-1)^{4-l} l! \genfrac{\{}{\}}{0pt}{}{4}{l} 
B_n^{\overbrace {\scriptstyle{(0,\cdots,0}}^{l+1})}+
\displaystyle \sum_{m=1}^{2} (-1)^{2-m} a_{1,m} B_n^{(0,0,-m)}\\
&=-\genfrac{\{}{\}}{0pt}{}{4}{1}B_n^{(0,0)}+2!\genfrac{\{}{\}}{0pt}{}{4}{2}B_n^{(0,0,0)}
-a_{1,1}B_n^{(0,0,-1)}+a_{1,2}B_n^{(0,0,-2)}\\
&=-B_n^{(0,0)}+14B_n^{(0,0,0)}-10B_n^{(0,0,-1)}+2B_n^{(0,0,-2)}.
\end{align*}

\vspace{-8pt}
\hspace{-5pt}($a_{1,1}=a_{0,0}+2\cdot a_{0,1}=6+2\cdot 2=10$)

Moreover, we see the parts of coefficients on Multi-Poly-Bernoulli numbers of the identity which hold on Theorem3.14.
Therefore we revisit Example3.15.

\textbf{Example3.15}(Example3.15 revisited). 

We give examples of Theorem3.14 for $1\leq r\leq2$,\,$1\leq k \leq 6$.

(i)For $r=1$

$B_n^{(-1)}=B_n^{(0,0)}$

$B_n^{(-2)}=-B_n^{(0,0)}+2B_n^{(-1,0)}$

$B_n^{(-3)}=B_n^{(0,0)}-3B_n^{(-1,0)}+3B_n^{(-2,0)}$

$B_n^{(-4)}=-B_n^{(0,0)}+4B_n^{(-1,0)}-6B_n^{(-2,0)}+4B_n^{(-3,0)}$

$B_n^{(-5)}=B_n^{(0,0)}-5B_n^{(-1,0)}+10B_n^{(-2,0)}-10B_n^{(-3,0)}+5B_n^{(-4,0)}$

$B_n^{(-6)}=-B_n^{(0,0)}+6B_n^{(-1,0)}-15B_n^{(-2,0)}+20B_n^{(-3,0)}-15B_n^{(-4,0)}+6B_n^{(-5,0)}$

(ii)For $r=2$

$B_n^{(-1)}=B_n^{(0,0)}$

$B_n^{(-2)}=-B_n^{(0,0)}+2B_n^{(0,0,0)}$

$B_n^{(-3)}=B_n^{(0,0)}-6B_n^{(0,0,0)}+6B_n^{(-1,0,0)}$

$B_n^{(-4)}=-B_n^{(0,0)}+14B_n^{(0,0,0)}-24B_n^{(-1,0,0)}+12B_n^{(-2,0,0)}$

$B_n^{(-5)}=B_n^{(0,0)}-30B_n^{(0,0,0)}+70B_n^{(-1,0,0)}-60B_n^{(-2,0,0)}+20B_n^{(-3,0,0)}$

$B_n^{(-6)}=-B_n^{(0,0)}+62B_n^{(0,0,0)}-180B_n^{(-1,0,0)}+210B_n^{(-2,0,0)}-120B_n^{(-3,0,0)}+30B_n^{(-4,0,0)}$

Here by Example3.15, the sum of coefficients on Multi-Poly-Bernoulli numbers are all $1$ for $r=1$ and $r=2$.
From this, we can be considered the following Theorem3.18.

\textbf{Theorem3.18.} We have the following relations on the sum of coefficients.

$\displaystyle \sum_{l=1}^{r} (-1)^{k-l} l! \genfrac{\{}{\}}{0pt}{}{k}{l}
+\displaystyle \sum_{m=1}^{k-r}
(-1)^{k-m-r}  \binom{k}{m} r! \genfrac{\{}{\}}{0pt}{}{k-m}{r}=1_{\textbf.}$

\begin{align*}
\mathrm{Proof.} \displaystyle \sum_{m=1}^{k-r}
(-1)^{k-m-r}  \binom{k}{m} r! \genfrac{\{}{\}}{0pt}{}{k-m}{r}
&=\displaystyle \sum_{m=r+1}^{k}
(-1)^{k-m}  \binom{k}{m-r} r! \genfrac{\{}{\}}{0pt}{}{k-m+r}{r}\\
&=\displaystyle \sum_{l=r+1}^{k}
(-1)^{k-l} l! \genfrac{\{}{\}}{0pt}{}{k}{l}.
\end{align*}
Thus, the left hand of the equality equals
$\displaystyle \sum_{l=1}^{k} (-1)^{k-l} l! \genfrac{\{}{\}}{0pt}{}{k}{l}$,
and from the proof of Theorem3.12, we have
$\displaystyle \sum_{l=1}^{k} (-1)^{k-l} l! \genfrac{\{}{\}}{0pt}{}{k}{l}$=$1$.

Therefore we obtain the result. \hspace{8cm} ${\Box}$

Similary, we see the parts of coefficients on Multi-Poly-Bernoulli numbers where represents $B_n^{(-k)}$ by using the sum of $B_n^{\overbrace {\scriptstyle{(0,\ldots,0}}^{l+1})}$ \hspace{3pt}
($ l \geq 1$) and them of $B_n^{\overbrace {\scriptstyle{(0,\ldots,0}}^{r},-m)}$. But we don't find regularities, therefore we don't see the relations yet.

We give tables which show the values of $B_n^{(k)}$  ($-5 \leq k\leq 5$,\hspace{3pt}$0 \leq n \leq 7$) and $B_n^{(k_1,k_2)}$
for small $n$, $k_i$.

\begin{center}
Table 2. [1] $B_n^{(k)}$ ($-5 \leq k\leq 5$,\hspace{3pt}$0 \leq n \leq 7$)

\begin{tabular}{|c|c|c|c|c|c|c|c|c|} \hline
\hspace{5pt}k $\backslash$ n&0&1&2&3&4&5&6&7 \\ \hline
-5 &1 &32 &454 &4718 &41506 &329462 &2441314 &17234438 \\ \hline
-4 &1 &16 &146 &1066 &6902 &41506 &237686 &1315666 \\ \hline
-3 &1 &8 &46 &230 &1066 &4718 &20266 &85310 \\ \hline
-2 &1 &4 &14 &46 &146 &454 &1394 &4246 \\ \hline
-1 &1 &2 &4 &8 &16 &32 &64 &128 \\ \hline
0 &1 &1 &1 &1 &1 &1 &1 &1 \\ \hline
1 &1 &$\frac{1}{2}$ &$\frac{1}{6}$ &0 &$-\frac{1}{30}$ &0 &$\frac{1}{42}$ &0 \\[5pt] \hline

2 &1 &$\frac{1}{4}$ &$-\frac{1}{36}$ &$-\frac{1}{24}$ &$\frac{7}{450}$ &$\frac{1}{40}$ &$-\frac{38}{2205}$ &$-\frac{5}{168}$ \\[5pt] \hline

3 &1 &$\frac{1}{8}$ &$-\frac{11}{216}$ &$-\frac{1}{288}$ &$\frac{1243}{54000}$ &$-\frac{49}{7200}$ &$-\frac{75613}{3704400}$ &$\frac{599}{35280}$ \\[5pt] \hline

4 &1 &$\frac{1}{16}$ &$-\frac{49}{1296}$ &$\frac{41}{3456}$ &$\frac{26291}{3240000}$ &$-\frac{1921}{144000}$ &$\frac{845233}{1555848000}$ &$\frac{1048349}{59270400}$ \\[5pt] \hline

5 &1 &$\frac{1}{32}$ &$-\frac{179}{7776}$ &$\frac{515}{41472}$ &$-\frac{216383}{194400000}$ &$-\frac{183781}{25920000}$ &$\frac{4644828199}{653456160000}$ &$\frac{153375307}{49787136000}$ \\[5pt] \hline
\end{tabular}
\end{center}

\begin{center}
Table3. [5] $B_n^{(k_1,k_2)}$ ($0 \leq n \leq 7$, $k_1, k_2$: small values)

\begin{tabular}{|c|c|c|c|c|c|c|c|c|} \hline
\hspace{5pt} $\backslash$ n&0&1&2&3&4&5&6&7\\ \hline
$B_n^{(1,1)}$&$\frac{1}{2}$&$\frac{1}{2}$&$\frac{5}{12}$&$\frac{1}{4}$&$\frac{1}{20}$&$-\frac{1}{12}$&$\frac{5}{84}$&$\frac{1}{12}$\\ \hline
$B_n^{(1,0)}$&1&$\frac{3}{2}$&$\frac{13}{6}$&3&$\frac{119}{30}$&5&$\frac{253}{42}$&7\\ \hline
$B_n^{(0,1)}$&$\frac{1}{2}$&$\frac{2}{3}$&$\frac{5}{6}$&$\frac{29}{30}$&$\frac{31}{30}$&$\frac{43}{42}$&$\frac{41}{42}$&$\frac{29}{30}$\\ \hline
$B_n^{(0,0)}$&1&2&4&8&16&32&64&128\\ \hline
$B_n^{(0,-1)}$&2&6&18&54&162&486&1458&4374\\ \hline
$B_n^{(-1,0)}$&1&3&9&27&81&243&729&2187\\ \hline
$B_n^{(-1,-1)}$&2&9&39&165&687&2829&11505&46965\\ \hline
\end{tabular}
\end{center}

\textbf{References}

[1]\, T.Arakawa, T.Ibukiyama\,and M.Kaneko, Bernoulli Numbers and Zeta 

\hspace{0.5cm} Functions, Springer, Tokyo, (2014)

[2]\, C. Brewbaker. A combinatorial interpretation of the poly-Bernoulli 

\hspace{0.5cm} numbers and two Fermat analogues. Integers 8\,(2008), A02.

[3]\, S. Launois. Rank t H-primes in quantum matrices. Comm. Algebra 

\hspace{0.5cm} 33\,(2005), 837-854.

[4]\, K.Kamano,\,A\,formula\,for\,Multi-Poly-Bernoulli\,numbers\,of\,negative\,index.

\hspace{0.5cm} Kyushu\,J. Math.67  (2013),\,29-37

[5]\, Y.Hamahata\,and\,H.Masubuchi,\,Recurrence\,formulae\,for\,Multi-Poly-Bern-

\hspace{0.5cm} oulli numbers. Integers\,7(2007),\,A46,\,1-15

[6]\, Graham,R.,Knuth,D.,Patashnik,O.:Concrete\,Mathematics.\,Addison-

\hspace{0.5cm} wesley(1989)

\end{document}